\documentclass[onefignum,onetabnum]{siamonline220329}

\usepackage{lipsum}
\usepackage{amsfonts}
\usepackage{graphicx}
\usepackage{epstopdf}
\ifpdf
  \DeclareGraphicsExtensions{.eps,.pdf,.png,.jpg}
\else
  \DeclareGraphicsExtensions{.eps}
\fi

\usepackage{enumitem}
\setlist[enumerate]{leftmargin=.5in}
\setlist[itemize]{leftmargin=.5in}

\newsiamremark{remark}{Remark}
\newsiamremark{hypothesis}{Hypothesis}
\crefname{hypothesis}{Hypothesis}{Hypotheses}
\newsiamthm{claim}{Claim}

\usepackage{amsopn}

\usepackage{amsmath,amssymb,epsfig,float,color}
\usepackage{tabularx} 
\usepackage[noend]{algpseudocode}   

\newcommand{\R}{\mathbb{R}}

\renewcommand{\refname}{}
\newcommand{\half}{{\frac{1}{2}}}

\newcommand\hl[1]{%
  \bgroup
  \hskip0pt\color{red!80!black}%
  #1%
  \egroup
}

\newcommand{\uhat}{\hat{u}}

\newcommand{\grad}{\nabla}

\def\O{{\it O}}

\def\ndata{\hat{n}}

\headers{Model discovery }{T. Meissner and K. Glasner}

\title{Simultaneous model discovery and state estimation under high data corruption}

\author{Teddy Meissner\thanks{Corresponding Author: Program in Applied Mathematics, University of Arizona, Tucson AZ 85721 USA (\email{tmeissner@arizona.edu}).}
\and Karl Glasner \thanks{Department of Mathematics, University of Arizona, Tucson AZ 85721 USA (\email{kglasner@arizona.edu}).} }

\begin{document}

\maketitle

\begin{abstract}
This paper proposes a sparse regression strategy for discovery of ordinary differential equations from incomplete and noisy data.  Inference is performed over both equation parameters and state variables
using a statistically motivated likelihood function.
Sparsity is enforced by a selection algorithm which iteratively removes terms and compares models using statistical information criteria.  Large scale optimization is performed using a second-order variant of the Levenberg-Marquardt method, where the gradient and Hessian are computed via automatic differentiation.   The proposed method is illustrated and tested on several systems with varying levels of noisy and incomplete data.   Comparisons are made to a state-of-the-art algorithm
for system identification, demonstrating competitiveness of the proposed approach.

\end{abstract}




\section{Introduction}

Combining physical laws with observational data is an ongoing scientific challenge.
Whereas differential equations derived entirely from first principles have been 
broadly successful in basic physics and mechanics, this approach is insufficient for complex systems ranging
from materials science to biology.  In those cases, model development involves a variety of
phenomenological assumptions, which must be calibrated or verified from observations.
Many recent efforts have focused on the problem of model discovery, which
involves generating equations from scratch that attempt to balance trade-offs between fitting data,
parsimony, and interpretability of the discovered model.  

Historically, the earliest methods of combining data and differential equations were
concerned with the inverse problem of parameter estimation (e.g. \cite{ramsay2005fitting} and references therein).  
Traditional nonlinear least squares regression (known in various contexts as the shooting method or adjoint method),
seeks to minimize only the difference between model output and data (e.g.  \cite{fullana1997identification,ackleh1998numerical,ashyraliyev2008parameter,garvie2010efficient,jin2012sparsity,croft2015parameter,ramsay2017dynamic,sgura2019parameter,zhao2020learning,glasner2021optimization}).
An alternative is to optimize over the equation's residual, the error obtained when the data is inserted into the model.  This idea forms the basis for many recent studies of both parameter estimation and model discovery  
\cite{ramsay2017dynamic,brunel2008parameter,gugushvili2012sqrt,liang2008parameter,brunton2016discovering,schaeffer2017learning}.
Methods which simultaneously infer the state and parameter variables (hereafter called ``hybrid" regression) interpolate between these two approaches 
\cite{ramsay2007parameter,voss2004nonlinear,galioto2020bayesian,retzler2022shooting,ribera2022model,glasner2023data}, and have proven to be more robust for handling noisy and incomplete data.
Gaussian process regression  \cite{raissi2018hidden,raissi2017machine} has also been proposed for system identification.  Lastly, a variety of studies employ Bayesian inference 
(e.g. \cite{dewar2010parameter,barajas2019approximate,galioto2020bayesian,yoshinaga2020bayesian,zhao2020learning,zhao2021image}) that can provide a detailed probability distribution of parameter values, which is useful for tasks such as uncertainty quantification.

Discovering differential equations from data, rather than estimation of parameters in
a known model, has seen a surge of interest. 
Some of the earliest work in this direction involves symbolic regression 
\cite{bongard2007automated,schmidt2009distilling}.
These procedures consist of generating candidate equations by symbolically combining elementary operations,
and testing them competitively to ensure a good fit with the data.  A variety of extensions build upon this idea, employing either genetic algorithms or deep learning strategies
\cite{long2019pde,udrescu2020ai,atkinson2019data,maslyaev2021partial,chen2022symbolic}.

An alternative to symbolic discovery involves regression with a given set of candidate terms, or at least a general equation structure. 
In order to identify models which are both parsimonious and do not overfit data, these methods rely on promoting sparsity
in the parameter set.  This can be accomplished by using statistical information criteria \cite{mangan2017model}, $L^1$ penalization
\cite{schaeffer2017learning}, or a sequential thresholding algorithm (``SINDy")
\cite{brunton2016discovering,rudy2017data}.  
Sparsity has also been incorporated into Bayesian versions of model discovery
\cite{bonneville2022bayesian,hirsh2022sparsifying}.

It has become increasingly popular to employ neural network function approximation
\cite{berg2017neural,raissi2019physics} instead of traditional grid-based or finite element representations.
These ``physics-informed" machine learning methods \cite{karniadakis2021physics} rely on a 
loss function similar to hybrid regression methods discussed above, inferring
both the solution and parameters.
They have been adapted to model discovery by including a sparsity-inducing penalty term in the loss function,
and search for equations composed of predefined library terms \cite{both2021deepmod,xu2021deep}
or represented by rational neural networks \cite{stephany2022pderead,stephany2022pdelearn}.
Alternatively, neural networks may be employed to approximate functions that define a dynamic
model directly \cite{kiyani2022machine,haring2022levenberg,nguyen2020assimilation}.

Despite the proliferation of methods for data-driven model discovery, there are a variety of outstanding challenges.  Here we highlight some of these and explain how they are addressed here. 
\begin{itemize}
\item
{\bf Noisy data.}
Methods that utilize residual regression, most notably the SINDy algorithm, are very sensitive to noise. Subsequent modifications of the original algorithm address this by
direct filtering of the data \cite{schaeffer2017sparse,wentz2023derivative}, 
ensemble techniques which average results using multiple data sets \cite{fasel2022ensemble,delahunt2022toolkit}, or employ additional sparsity requirements
in the data set so as to identify outliers \cite{zhang2023adaptive}.  
Another body of work utilizes  Galerkin-type weak formulations of the equations
\cite{reinbold2020using,messenger2021weak,messenger2021weak2,tang2023weakident},
which produces the effect of a customized low-pass filter.  
Our method does not rely on direct or indirect filtering, but rather uses inference to produce
a system state which interpolates between the provided data and an exact smooth solution of the differential equation, which is determined simultaneously.

\item
{\bf Incomplete data.}     
For practical reasons, observed output values of the model may not be present throughout the domain, or might not be provided for every solution component.  
This is clearly problematic for residual regression approaches which rely on differentiation of the data itself, although one proposed weak formulation shows
promise in handling latent (unmeasured) variables  \cite{reinbold2020using}.
The hybrid regression approach used here, on the other hand, looks to fit a model with whatever data is available, and can be effective for inference of differential equations
using remarkably little information \cite{ramsay2007parameter,ribera2022model,glasner2023data}.

\item
{\bf State inference.}  
In the course of model discovery, one might also want to estimate the true ground state, that is the solution of the actual model.  This can in principle be done with any discovery method for which
side conditions are known, simply by integrating the discovered system forward in time.
This would be problematic in chaotic systems, however, since even if the exact model were known, solution trajectories are highly sensitive to noise in the
initial condition and would diverge from the true state.
Our method simultaneously infers both the model and a state which represents an approximate solution to the discovered equations.
Our numerical experiments (section \ref{sec:Lorenz})  verify that the ground state solution can be accurately recovered in chaotic systems even in the presence of significant noise.
 
\item
{\bf Equations which are nonlinear in parameters.}
A typical staring point for model discovery involves choosing a library of candidate terms,
and supposing the equations are sparse linear combinations of these terms.   This is necessary for methods that rely on residual regression (e.g. SINDy and its variants), since the least squares problem is exactly linear.  On the other hand, this assumption is fairly restrictive, and a more general regression method is needed to deal with such cases (see section \ref{Sec:Colpitts}).
\end{itemize}

We employ a formulation that accounts for uncertainties in both the model and data,
an idea which has statistical roots in data assimilation 
\cite{evensen279data,ye2015improved}.
In this framework, ignoring uncertainties in the data leads to (SINDy-type) residual regression \cite{galioto2020bayesian}, which explains why such approaches struggle
with noisy data.  Conversely, parameter estimation that ignores model uncertainty can be highly non-convex and problematic from an optimization standpoint \cite{ramsay2005fitting,glasner2023data}.

Recently, system identification methods have been proposed which take both types of
uncertainly into account.   This is most commonly done in conjunction with
maximum likelihood estimation \cite{ramsay2007parameter,voss2004nonlinear,ribera2022model,glasner2023data}, and results in a hybrid regression method that optimizes over 
both parameters and state variables.  The multiple shooting method of 
Retzler et al. \cite{retzler2022shooting} is a version of this, 
where only a small subset of state variables are inferred, and the remainder are determined by forward integration.  A more elaborate Bayesian setting was utilized
by Galioto \& Gorodetsky \cite{galioto2020bayesian}, who use an exact likelihood expression that accounts for both types of uncertainty, but 
leads to computationally expensive integrals that ultimately require Monte Carlo techniques.   The maximum likelihood approach was extended to include 
model sparsity by Ribera et al. \cite{ribera2022model}.
Lastly, sparse system identification using neural networks has been 
proposed by Nguyen et al. \cite{nguyen2020assimilation} and Stephany and Earls \cite{stephany2022pderead,stephany2022pdelearn}.  These works employ a loss function
similar in character to the objective in maximum likelihood estimation.

This paper introduces a hybrid sparse regression approach that accounts for data and model uncertainty.  It is shown to accurately reconstruct systems of ordinary differential equations in circumstances where the supplied data is both noisy
and incomplete.  Novel aspects include:
\begin{itemize}
\item
Sparsity is enforced via Bayesian information criteria, rather than by sweeping over hard threshold values as in \cite{ribera2022model};
\item
We demonstrate that a second-order optimization method can be made practical using
automatic differentiation and a graph-coloring algorithm that takes advantage of the sparsity of the Hessian matrix.  This contrasts to first order optimization methods used elsewhere for system identification
\cite{ribera2022model,nguyen2020assimilation,stephany2022pderead,stephany2022pdelearn}.
\item
Extensive numerical experiments are conducted, investigating the proposed method's
robustness to high noise and missing data.
\item
Comparisons of both system recovery rates and parameter estimates are made to the WSINDy algorithm \cite{messenger2021weak},  a widely used benchmark
for system identification.
\item
Our approach is formulated with generality in mind.  In particular,
there is no restriction to equations which are linear in parameters, which is a common drawback to many library-based methods.  We provide an example in section \ref{Sec:Colpitts} to illustrate this feature.
\end{itemize}

This paper is structured as follows.  Formulation of our regression strategy and its statistical underpinnings are given in section \ref{sec:setup}.   The sparse selection methodology that leads to an explicit algorithm is discussed in section \ref{sec:methodology}.  
Details of numerical discretization and optimization methods employed are provided in section \ref{sec:optimization}.  Finally, several computational examples and comparisons of model discovery
are given in section \ref{sec:results}.

\section{Problem formulation}
   \label{sec:setup}
Suppose that a differential equation is discretized and written in the form
\begin{equation}
    N(u;\theta) = 0,
\end{equation}
where $u \in \R^n$ is the state vector and $\theta \in \R^p$ encapsulates all the unknown constitutive parameters in the model, which might arise within the equation or side conditions. The function $N: \R^n \times \R^p \to \R^n$ is assumed to be a generic, differentiable nonlinear operator. No particular type of discretization (finite difference, finite element, neural network, etc.) is assumed at this stage, only that the approximate solution is described by the vector $u$.

Our main interest is to determine a sparse estimate for $\theta$ inferred from information about the solution, which could be obtained either from physical observation or derived from other computations.  This information is supplied in the form of noisy and possibly incomplete data $\hat{u} \in \mathbb{R}^{\hat{n}}$, which approximates $u$ on a subset $D$ of components.  The notation $v|_D$ will represent restriction of the vector $v$ to the subset $D$.

We begin with the simpler problem of estimation of $\theta$ without any regard to
sparse selection. A classical approach to parameter estimation is maximum likelihood regression, which utilizes a likelihood function determined from a known (or parametric) statistical distribution of data.
For example, if the only source of uncertainty is zero-mean Gaussian noise with variance $\sigma^2$ in the data and $u^*(\theta)$ represents the exact (and unique) solution for parameters $\theta$, the likelihood is the conditional probability
\begin{equation}
  \label{leastsq}
   L(\theta) = (\sigma\sqrt{2\pi})^{-\hat{n}}
   \exp \left(-\frac{\|\hat{u}-u^*(\theta)|_D\|_2^2}{2\sigma^2} \right).
\end{equation}
Estimation of $\theta$ is therefore a nonlinear least squares fitting procedure,
obtained by maximizing (\ref{leastsq}). 

A more general approach takes into account uncertainty in both the data and the underlying model.  For instance,
$N(u;\theta)$ can be regarded as random, obeying a distribution with zero-mean and small variance. 
Under these assumptions,  an approximate log-likelihood function can be written as \cite{glasner2023data}
\begin{equation}
  \label{like2}
L(\theta) = - \frac{\ndata}{2} \ln F(\theta), \quad
F(\theta) \equiv  \min_u  \lambda || \, \uhat - u|_D  \, ||_2^2 + ||N(u;\theta)||_2^2.
\end{equation}
Here $\lambda$ is a ratio of variances that describes the relative uncertainty in the model
compared to the data. Expression (\ref{like2}) infers both parameter and solution variables
$\theta$ and $u$ in the following sense.
Defining
\begin{equation}\label{LossNoSparse}
    \mathcal{L}_{\lambda}(\theta,u) = \|N(u;\theta)\|_2^2 + \lambda \| \, \hat{u} - u|_D \, \|_2^2,
\end{equation}
the maximum likelihood estimate is obtained by minimizing $\mathcal{L}_{\lambda}$ over
both $\theta$ and $u$.  

A brief study of the limiting behavior in (\ref{LossNoSparse}) for small and large $\lambda$ 
is illuminating, as we can see that this optimization problem interpolates between two types of regression:  nonlinear least squares, and residual-only minimization.
For $\lambda \ll 1$, the leading order minimization problem
involves only the first term in  (\ref{LossNoSparse}).  This problem is degenerate,
and has a family of global minimizers $u^*(\theta)$ which solve  $N(u^*(\theta),\theta)=0$.
By expanding the exact minimum $u_{min} = u^*(\theta) + \O(\lambda)$, a solvability condition is determined
at order $\lambda$, which involves minimization of 
$  \| \, \hat{u} - u^*(\theta) |_D \, \|_2^2$ over $\theta$.  This problem is precisely the same as minimizing (\ref{leastsq}), which is the nonlinear least squares solution sought in classical parameter estimation.

In contrast, when $\lambda \gg 1$, the second term in (\ref{LossNoSparse}) dominates,
and the leading order state variable is $u \approx \hat{u}$ (assuming no restriction to
a subset $D$ of data).
This is again a degenerate situation, as $\theta$ is not determined at this level of optimization.  
In this case, the expansion $u_{min} = \hat{u} + \lambda^{-1} u_1 + \ldots$ leads to the correction defined by     
\begin{equation}
   \label{O2}
  \min_{u_1,\theta}   \|N(\hat{u}; \theta)\|_2^2 +  \| \, \hat{u} - u_1 \, \|_2^2.
\end{equation}
This problem decouples, leading to the residual minimization problem
\begin{equation}  \min_{\theta}  \|N(\hat{u}; \theta)\|_2^2. 
\end{equation}
This type of regression is at the heart of the SINDy algorithm \cite{brunton2016discovering} and its
derivatives.  From a statistical point of view, this limit ignores uncertainly in the data, which is consistent with the observation that residual regression struggles with noise.  Various attempts to improve this situation have been proposed, significantly the weak-SINDy algorithm \cite{messenger2021weak}, which we later make empirical comparisons to.

Discovery of equations, as opposed to parameter estimation, must take into account a wide variety of candidate models as well as parameters.  Model selection can be performed by insisting on sparsity, that is, selecting a small number of terms which provides a reasonable description of the supplied data.
A typical approach is to include a sparsity promoting penalty $R(\theta)$ in the regression formulation, 
 \begin{equation}\label{LossSparse}
    \mathcal{L}_{\lambda,R}(\theta,u) = \|N(u;\theta)\|_2^2 + \lambda \|\hat{u} - u|_D\|_2^2 + R(\theta).
\end{equation}
Various choices for the penalty function are described in the next section.

\section{Methodology} \label{sec:methodology}
The problem setup described in section \ref{sec:setup} makes no restriction on the choice of 
generalized model $N(u;\theta)$.  Commonly, $N(\cdot)$ is taken to be a linear combination
\begin{equation}
  \label{lincomb}
  N = N_0(u) + \sum_{k=1}^p \theta_k  N_k(u).
\end{equation}
where $\{ N_k \}$ represent a predetermined library of possible candidates.
We will also consider a case (section \ref{Sec:Colpitts}) where the library terms themselves have parameter dependence.

A constraining feature of library-based system identification methods (SINDy, etc.)
is the need for subjective assumptions about the structure of the system.  This could be
done using expert knowledge of a scientific application, for example by tailoring them to obey symmetries or underlying physics.  More commonly, library terms are specified generically, e.g. as multinomials, which is the approach taken in our numerical studies.  While in principle multinomials can approximate functions in a general way, in the spirit of parsimonious modeling it is reasonable to limit their order to a reasonable predefined size.  

It is supposed that most terms in (\ref{lincomb}) have no significance in describing the data.
One approach to identifying such terms is the addition of a sparsity-inducing penalty
as in (\ref{LossSparse}).    A common choice for $R(\theta)$ is provided by the ``norm" $||\theta||_0$, which is a simple count of the number of non-zero components.
Despite it being the most natural measure of sparsity, it is problematic in optimization, 
and can lead to intractably high computational complexity. 
Alternatives include other norms on $\theta$, notably
$l_1$ or $l_2$, which do not exhibit the same computational
difficulties associated with the $l_0$ penalty,
but are known to produce significant bias in parameter estimates.  

In our work, we utilize (\ref{LossSparse}) as a means to identify candidate terms which 
may or may not be removed according to an information criteria-based algorithm
(see section \ref{sec:optimizemodelAlgorithm}).  We adopt the regularized $l_0$ penalty
proposed in Mohimani et al. \cite{mohimani_fast_2009}
\begin{equation}
\|\theta\|_\epsilon = \sum_{i=1}^d \left(1 - \exp\left(-\frac{\theta_i^2}{2\epsilon^2}\right)\right),
\end{equation}
where it is clear that $\lim_{\epsilon\rightarrow 0}\| \theta \|_\epsilon = \| \theta \|_0$.  This smooth approximation allows for gradient-based optimization methods for (\ref{LossSparse}), enabling more practical implementations in high-dimensional settings.

\subsection{Model comparison and subset selection}
\label{sec:selection}
Under assumptions of model sparsity,  it is necessary to prescribe methods for selecting and comparing different subsets of model terms.  The most straightforward approach is to 
recursively remove terms with parameters below a predefined threshold value
(e.g. \cite{brunton2016discovering,rudy2017data}).  Alternatively, accepting the removal of one or more terms may be decided by a different measure.  Two well-established model selection tools used in statistical analysis and machine learning are the Akaike Information Criteria (AIC) \cite{Akaike1998-zv} and the closely-related Bayesian information criteria (BIC) \cite{Schwarz1978-qx}.  Both of these seek to balance the goodness of fit of the model with its complexity.  The practical difference is that BIC places a higher penalty on models with more parameters, and it is more stringent against overfitting in scenarios with large datasets. 
We found that empirically,  BIC gave somewhat more satisfactory results, and is used in our latter tests.

The Bayesian Information Criterion (BIC)  (also known as the Schwarz criterion)
is defined as
\begin{equation}
\text{BIC} = \ln(n)d - 2 \log L,
\end{equation}
where $d$ is the number of parameters in the model and $L$ is the likelihood of the model given the data.   Models with smaller BIC scores are preferred.
For our specific likelihood function (\ref{like2}), we therefore seek to minimize
\begin{equation} \label{BIC}
    \text{BIC}(u,\theta) = \ln(\hat{n})||\theta||_{0} + \hat{n}\ln \Big(
    \|N(u;\theta)\|_2^2 + \lambda \|\hat{u} - u|_D\|_2^2\Big)
\end{equation}
over both $u$ and $\theta$.   In practice, the second term is minimized for models
with a specified combination of terms, and the corresponding information criteria scores
are used to decide which is best.

\subsection{Model selection algorithm} \label{sec:optimizemodelAlgorithm}
Given a library with $p$ potential terms, there are $2^p$ subset models, so it is impractical to evaluate the information criteria of the entire set of sparse models.
Instead, our aim is to efficiently assess a limited range of possibilities, by selectively removing terms guided by the sparse loss function (\ref{LossSparse}),
and verifying that (\ref{BIC}) is reduced at each step.   An outline of the algorithm is presented in Algorithm \ref{Alg:StepwiseSelection}.

Our method begins with simultaneous minimization of the non-sparse loss \(\mathcal{L}_\lambda\) and the sparse loss \(\mathcal{L}_{\lambda, R}\). 
Initial guesses are provided for the first iteration of optimization, and thereafter are provided by previous optimizations.
The dual loss functions serve distinct purposes:  the unpenalized objective  \(\mathcal{L}_\lambda\) provides the likelihood of a model given a specific set of parameters and is used in (\ref{BIC}), whereas  \(\mathcal{L}_{\lambda, R}\) encourages irrelevant terms to be small, so they can be targeted for removal.  

\begin{figure}
\centering
\includegraphics[height = 6cm]{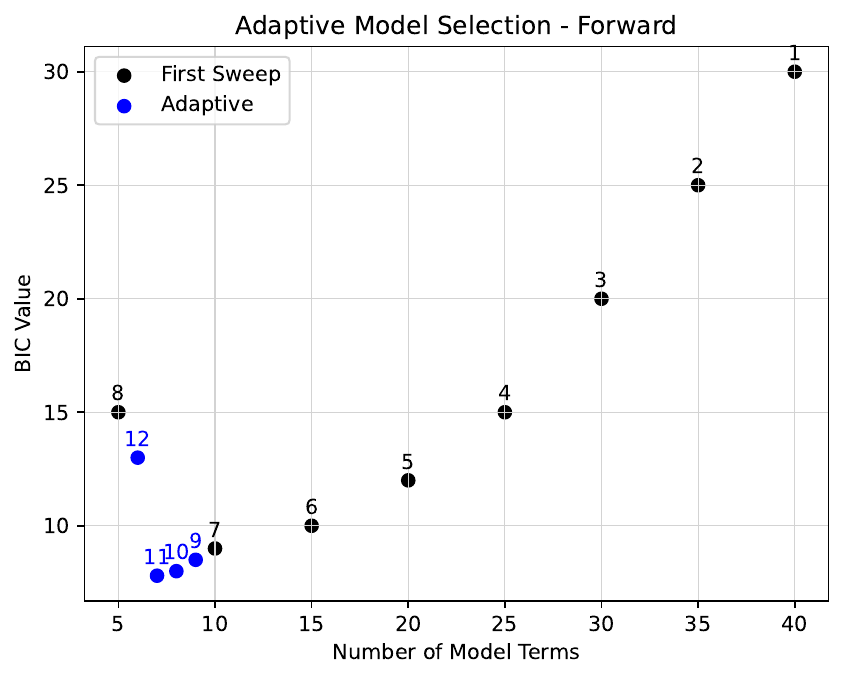}
\includegraphics[height = 6cm]{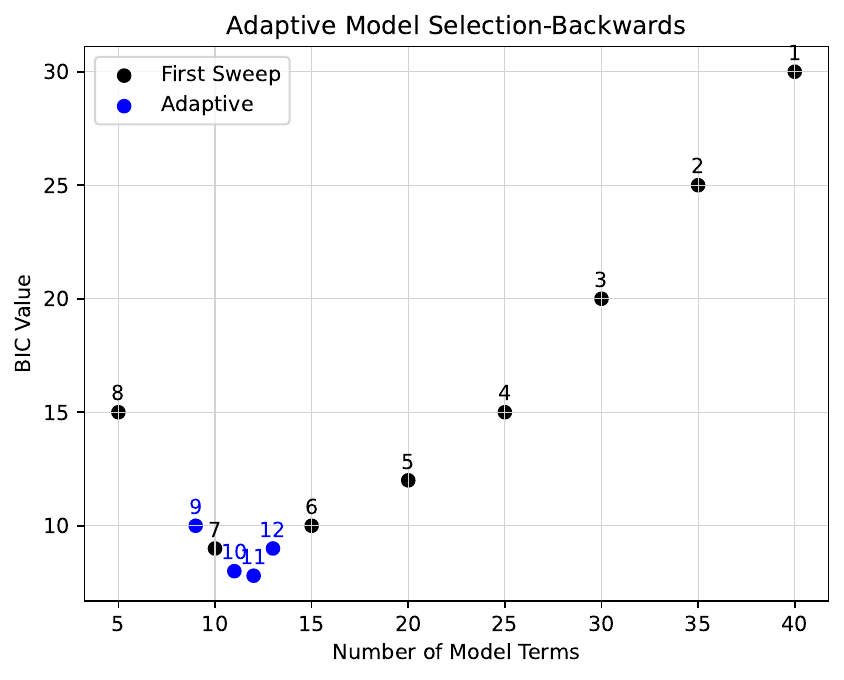}
\caption{Visualization of subset selection methodology to hone in on the correct terms under two different scenarios with parameter drop count $k = 5$. Iteration numbers are shown above each point, while both cases show the first model is rejected in iteration 8. Left shows acceptance of small forward steps resulting in no need to go to more complex models. Right shows a small forward step is rejected, resulting in the need for backward steps or testing more complex models.}
\label{fig:AICToy}
\end{figure}

The algorithm begins with all terms in the model library,  and selects an initial number $k=k_0$ of terms which will be removed at each iteration, provided that (\ref{BIC}) decreases.  If this is not the case,  the algorithm reverts to the best model so far and reduces $k$ to one for small refinements around the current best model. 
If further steps are accepted, the process continues in the {\it forward} refinement stage until (\ref{BIC}) increases, at which point the algorithm terminates
(see figure \ref{fig:AICToy}, left).
Alternatively, if pruning a single parameter from the best model leads to rejection just after $k$ is reduced,  the algorithm assumes that too many terms in the previous step were removed. 
This triggers the \emph{backwards} refinement stage  (see figure \ref{fig:AICToy}, right).  Here, the algorithm iteratively adds one term back into the model at a time, until the information criteria no longer decreases.  To determine which terms to add back in,  we revert to the  model which has $k_0$ more terms than the best model, that is, the iteration just prior
to the current minimum. From here, pruning $k_0 - 1$ terms is equivalent to adding one term back into the current model. The process continues by adding one term back at a time until no improvements are found. The algorithm requires at most $\lceil \frac{p}{k_0} \rceil + k_0 + 1$ iterations, where $\lceil \cdot \rceil$ denotes the ceiling function.

We can provide a heuristic argument about the convergence of the algorithm to the true minimum of (\ref{BIC}).   Define
$$
B(n)  = \min_{u,\theta, ||\theta||_0 = n} BIC(u,\theta),  
$$
where $n$ is less than or equal to the maximum number of library terms.   It follows that the minimum we seek is simply $\min_n B(n)$.
At each step of the algorithm,  we have a model with $n$ terms, and a corresponding value of (\ref{BIC}) we can call $\tilde{B}(n)$.  Provided
$\tilde{B}(n)$ is an accurate reflection of $B(n)$, minimizing $\tilde{B}(n)$ over $n$ would be sufficient.  Note that whether the algorithm finishes in the forward or backward stage,
the termination conditions are $\tilde{B}(n-1) \ge \tilde{B}(n^*)$ and $\tilde{B}(n+1) \ge \tilde{B}(n^*)$,  
assuming the final number of terms $n^*$ is strictly less than the library size.   If $B(n)$ is well approximated by $\tilde{B}(n)$, then the termination conditions guarantee that the desired minimum is recovered.

This approach is designed not only to identify the most parsimonious model, but also to methodically explore and document the performance of various models. Each tested model, along with its corresponding information criterion scores, is recorded, allowing for a comprehensive
analysis.  This strategy ensures that the search space is efficiently explored, increasing the likelihood of identifying the most effective model configuration.

\begin{algorithm}
\scriptsize
\caption{Adaptive Pruning with Forward Refinement and Backward Recovery}
\label{Alg:StepwiseSelection}
\begin{algorithmic}[1]
\Require Initial model $(u_0, \theta_0)$, block size $k$, data and sparsity weights $(\lambda, R)$
\State mask $M' = \mathbf{1}$ \Comment{Binary Mask}
\State mode $\gets$ forward, Adaptive $\gets$ True \Comment{Toggle variables for stepwise selection}
\State $(u^\star, \theta^\star,M^\star, IC^\star ) \gets (u_0, \theta_0, M', \infty)$ 
\State
\While{True}
    \State \textbf{in parallel:}
    \State \quad Solve $(u', \theta') \gets \arg\min \mathcal{L}_{\lambda}(u, \theta \odot M')$
    \State \quad Solve $(u^\dagger, \theta^\dagger)\gets \arg\min \mathcal{L}_{\lambda, R}(u, \theta \odot M')$
    \State $IC' \gets$ BIC$(u', \theta')$
    \If{$IC' < IC^\star$}
        \State Accept: $( u^\star, \theta^\star,M^*, IC^\star) \gets (M', u', \theta', IC')$
        \If {(Adaptive=True)$\land$(k = 1)} Adaptive $ \gets $ False \Comment{No backward refinement needed}  
        \EndIf
    \Else
        \If{(Adaptive = True)$\land$($k > 1$)}
            \State Reduce step: $k \gets 1$  \Comment{Try fine forward pruning}
        \ElsIf{(Adaptive=True)$\land$($k = 1$)}
            \State mode $\gets$ backward \Comment{Try fine backward pruning}
            \State Adaptive $ \gets $ False
        \Else
            \State \textbf{break}
        \EndIf
    \EndIf
    \State $M' \gets \textsc{UpdateMask}(\theta^\dagger, M^*, k, \text{mode})$
\EndWhile
\State \Return $(u^\star, \theta^\star, M^*)$
\\
\Function{UpdateMask}{$\theta^\dagger$, $M^\star$, $k$, mode}
    \If{mode = forward}
        \State Remove $k$ smallest terms in $\theta^\dagger$ not already masked
    \Else
        \State Add back $k$ largest terms previously masked
    \EndIf
    \State \Return new mask $M'$
\EndFunction
\end{algorithmic}
\end{algorithm}

\subsection{Hyperparameter selection}\label{sec:Validation}
The optimization problem in Algorithm \eqref{Alg:StepwiseSelection} depends on three parameters:  
\begin{itemize}
    \item $k_0$, which defines how many terms to prune each iteration,
    \item $\lambda$, which weights the fidelity to observed data, and  
    \item $R$, which promotes sparsity in the parameter vector $\theta$,
\end{itemize}
In this work we set $k_0 = 5$ for simplicity. We found this value to work well on a wide variety of problems, though we note tuning $k_0$ empirically had very little influence on the discovered model. Rather, it can be tuned based on the size of the library to influence fewer pruning updates.

\begin{figure}
    \centering
    \includegraphics[width=0.45\linewidth]{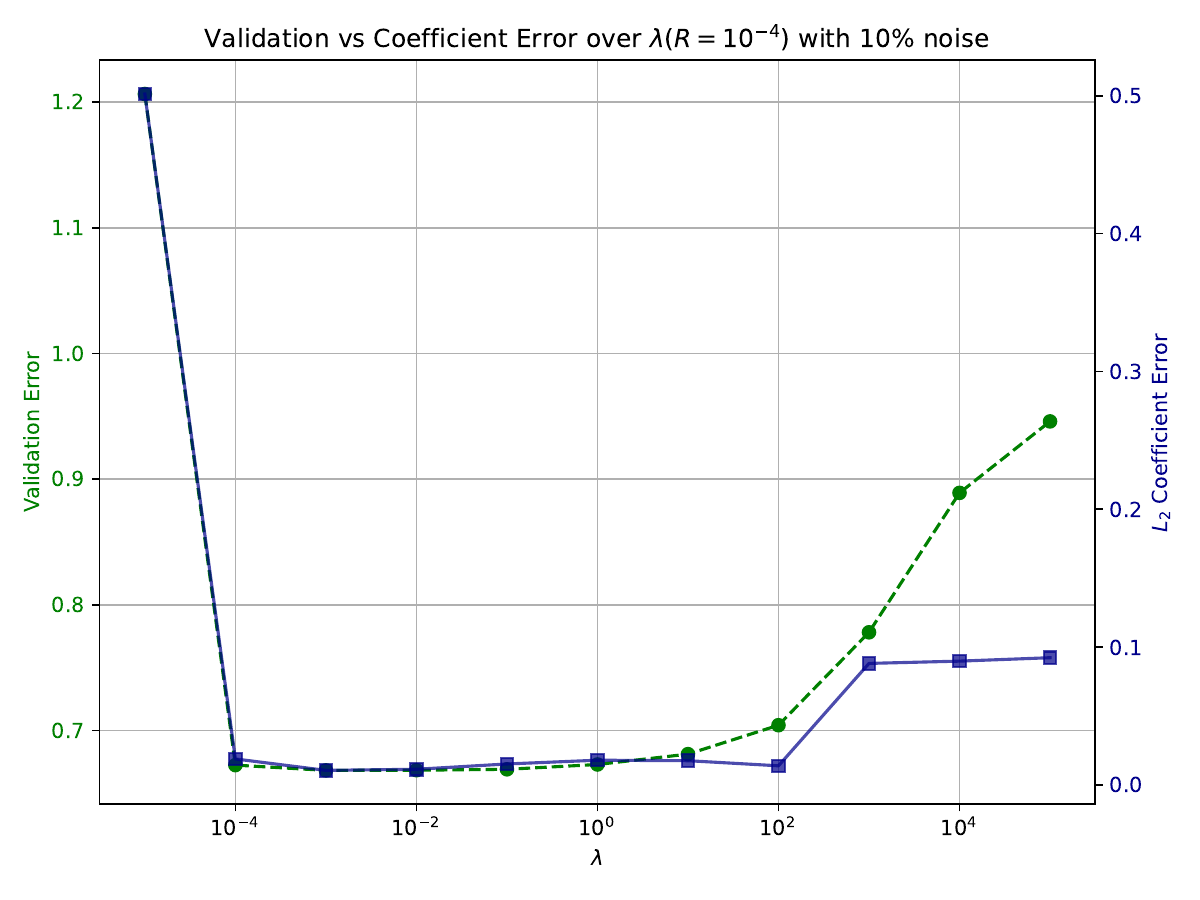}
    \includegraphics[width=0.45\linewidth]{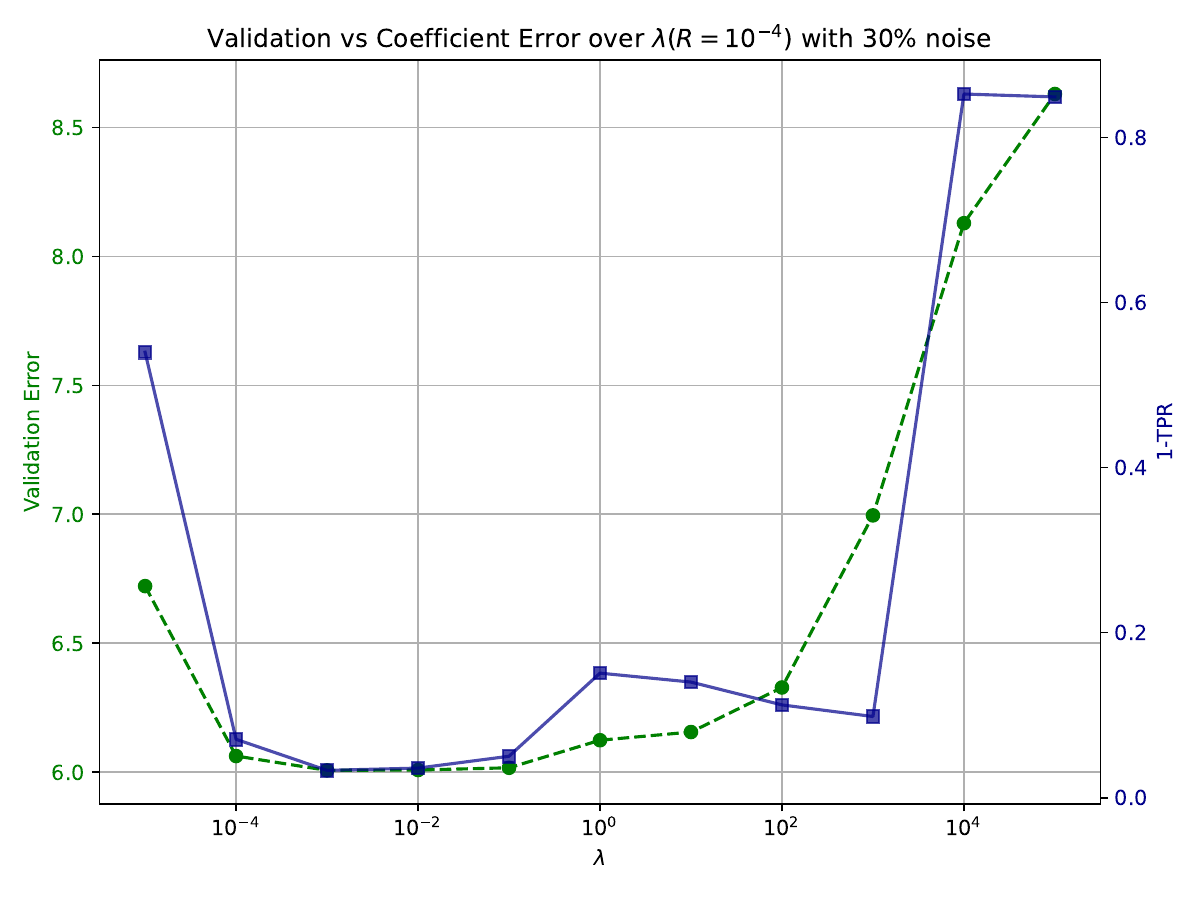}
    \caption{Validation error and relative $\ell_2$ error in the recovered coefficients for the Lorenz system under 10\% noise (left) and 30\% noise (right).}
    \label{fig:Cross-validation}
\end{figure}

To determine the best $(\lambda, R)$ pair, we use a cross-validation approach to find the weights which minimize the generalization error.
We split the data into two sets $D$ and $V$, where the former is used for optimization,
and the latter for validation.   Selection of $\lambda$ and $R$ is determined by 
\begin{equation}
    \min_{\lambda, R} \sum_{i \in V}  (\hat{u}_i - u_i)^2 
    \quad \text{s.t. $u$ solves Alg. \ref{Alg:StepwiseSelection} on $D$}
\end{equation}
In practice, we find that a wide range of hyperparameters prove satisfactory. We therefore enlist a simple grid search over a discrete set of candidate values.  For each pair $(\lambda, R)$,  the inner optimization problem is solved using Algorithm \ref{Alg:StepwiseSelection}, and the resulting model is evaluated with a fixed validation set. Figure \ref{fig:Cross-validation} presents a typical sweep over $\lambda$ values with $R$ held fixed, showing the validation error and relative $\ell_2$ error in the recovered coefficients compared to the ground truth, across varying noise levels. These results are based on the Lorenz system described in section \ref{sec:Lorenz}. While we display only $\lambda$ dependence in the figure, our full analysis tests combinations of both $\lambda$ and $R$. We see that the validation error and coefficient error tend to reach their minima at similar $\lambda$ values. This agreement indicates that minimizing validation error is an effective strategy for selecting models that not only generalize well but also recover accurate dynamics.

\subsection{Implementation}
\label{sec:implementation}

We are interested in discovery of ordinary differential systems
\begin{equation}
N(u; \theta) = \frac{du}{dt} - f(u; \theta) = 0, 
\end{equation}
where $u(t):\R \to \R^d$.   In the context of the sparse regression formulation,  
$N_0(u) \equiv \frac{d}{dt}u$ and 
\begin{equation}
  \label{ode}
   N(u; \theta) = N_0(u) - \sum_{k=1}^p \theta_k  N_k(u),
\end{equation}
where the linear combination of library terms has been identified with $f(u;\theta)$.
In practice, the state variable $u(t)$ is represented by some finite dimensional approximation,
and the time derivative is replaced with a numerical discretization formula.
It is also supposed that data $\{\hat{u}_i\}$ is provided at times $\hat{t}_1,\hat{t}_2,\ldots,\hat{t}_n$.  

The simplest way to represent the state variable is to consider values $u_i = u(t_i)$ only at sampling times 
$\{\hat{t}_i\}$.  On the other hand, 
if a more accurate inference of the state is required, one could either determine the state
from integration of the inferred system,  or a finer grid for values of $u(t)$ can be employed.
The former case leads to integration of an approximate system with approximate initial data, which can be wildly unstable;
this problem has been well studied for multiple shooting methods \cite{bock1983recent}.
In the latter case,  determining the appropriate grid size for $u(t)$ might require refinement studies,
or analytic error bounds, if they are available.    For comparisons of our approach to WSINDy
and investigations of incomplete data,  we suppose that the given data is sampled (aside from missing regions) at an acceptably fine scale which is used to construct the computational grid.  We alternatively study the effect of different computational grid sizes in relation to sampling rates (section \ref{subsection:Model_Refinement}).

Given a computational grid $\{t_i\}$, a general class of discretizations of the system (\ref{ode}) can be written
$$
N_{i+\half}  \equiv \frac{u_{i+1}  - I(u_i, \theta, \Delta t_i)}{\Delta t_i} = 0, \quad  \Delta t_i \equiv t_{i+1}-t_i,
$$
where $I$ is an integrator of choice using $u_i$ as the initial condition. We have found that the second-order midpoint method
\begin{align}
    \label{trapezoid}
    N_{i+\half}(\theta) &= \frac{u_{i+1} - u_{i}}{t_{i+1} - t_i} - 
    f \Big( \frac{u_{i+1} + u_{i}}{2}, \theta \Big) = 0, \quad i=1,\ldots,n-1,
\end{align}
provides a good trade-off between accuracy and simplicity.   

The data error term in
(\ref{LossNoSparse}) requires that $\{t_i\}$ includes a subset of values
$\{t_i | i \in D\}$ where the data is measured.  This is accomplished here by simply including the
sampling times $D$ in the discretized grid.  A reasonable alternative which is not investigated
in this work is to interpolate state values onto the set D.

With the discretization described above, the loss function $\mathcal{L}_{\lambda,R}(\theta,u)$ reads
\begin{equation}
   \label{loss_discrete}
    \mathcal{L}_{\lambda,R}(\theta,u) = \underbrace{\frac{1}{n}\sum_{i=1}^{n-1} N_{i+\half}(\theta)^2}_{\text{Model Error}} + \underbrace{\frac{\lambda}{\hat{n}} \sum_{i \in D} (\hat{u}_i - u_i)^2}_{\text{Data Error}} +  \underbrace{\frac{R}{k}\sum_{i=1}^k \left(1 - \exp\left(-\frac{\theta_k^2}{2\epsilon^2}\right)\right)}_{\text{Sparsity Penalty}},
\end{equation}
where $D$ is the set of indices in the domain at which data is available.   
Components of the ODE system can be written as $\{x_1, x_2, \dots, x_n\}$, and a typical
model library might look like
\begin{equation}
\label{ODE-Discrete}
N = \begin{bmatrix}
\dot{x}_1 \\
\dot{x}_2 \\
\vdots \\
\dot{x}_n
\end{bmatrix} - 
\begin{bmatrix}
x_1 & x_2 & x_1x_2 & x_1x_n & \cdots \\
x_1 & x_2 & x_1x_2 & x_1x_n & \cdots \\
\vdots & \vdots & \vdots & \vdots \\
x_1 & x_2 & x_1x_2 & x_1x_n & \cdots \\
\end{bmatrix}
\bigcirc 
\theta,
\end{equation}
where $M \bigcirc \theta$ denotes the element-wise Hadamard product of
a matrix $M$ and parameter matrix $\theta$, followed by summation over rows.

\section{Optimization}
            \label{sec:optimization}
The loss functions described in equations (\ref{LossNoSparse}) and (\ref{LossSparse}) are high-dimensional and generally non-convex, posing significant challenges to
optimization.  A wide variety of methods were investigated here, and we found
empirically that a generalized version of the Levenberg-Marquardt (LM) method using 
the exact Hessian had good performance (see section \ref{sec:comp}).

The original LM algorithm  \cite{Levenberg1944-qi,Marquardt1963-et} was developed strictly for least-squares problems,
and approximates the Hessian using a product of Jacobian matrices akin to the Gauss-Newton method.  We employ a more general version that uses the true Hessian \cite{Nocedal_undated-wf}, which essentially interpolates between gradient descent and Newton's method.
The precise details used in this paper are given in Appendix \ref{Appx: LM_opt}. 

The main computational effort here involves constructing the Hessian $H(x)$ and solving for the step $d_k$ by $H_\alpha d_k = - g$, where $H_\alpha = H + \alpha I$ is the modified Hessian and $g$ is the gradient. The parameter $\alpha$ is adjusted adaptively
to ensure that a local minimum is found.
In the proposed framework, the true Hessian is very sparse, allowing for efficient construction, detailed in section \ref{exploiting sparsity}. The time complexity is therefore dominated by the linear algebra.
In this work, the linear solve for $d_k$ uses sparse Cholesky factorization, while we note for very large problems
the conjugate gradient (CG) method may become a more suitable method.

\subsection{Computing the sparse Hessian}\label{exploiting sparsity}
In high-dimensional settings, efficient computation of the objective's Hessian matrix is a fundamental challenge.  Our approach benefits from sparsity that stems from the model's discretization, however, which makes second-order optimization feasible.  The loss functions are compatible with automatic differentiation, which simplifies construction of the Hessian across various applications. To effectively integrate sparsity and automatic differentiation, we first identify the sparsity structure of the Hessian, and subsequently apply graph coloring techniques to minimize the number of Hessian-vector products (HVPs) needed to compute the entire matrix.

\subsubsection{Finding the matrix structure}
To identify the sparsity pattern of the Hessian, we analyze the loss function in its discrete form given in \eqref{loss_discrete}. Considering each of the three components as independent functions, it is clear the data error and sparse penalty terms have no dependence on neighbor interactions, resulting in a contribution to the Hessian matrix only along the diagonal.  Derivatives of the form
 \begin{align}
     \frac{\partial^2 \|N(u;\theta)\|_2^2}{\partial u_m \partial \theta_n }      
\end{align}
are generally all non-zero, which creates a dense but small portion of the Hessian that is computed directly.

The largest portion of the Hessian corresponds to state derivatives,
\begin{align}
  \label{hess}
     \frac{\partial^2 \|N(u;\theta)\|_2^2}{\partial u_m \partial u_n } &=  \sum_{i,j \in N(m) \cap N(n)}  \frac{\partial^2 N(u_i; \theta)^2}{\partial u_i \partial u_{j} }.
\end{align}
where $N(x)$ denotes the neighbors connected at the $x^{th}$ grid point.  This relationship is provided by the finite difference stencils, which are known in advance and only computed once.  The neighbor relation can also be associated with a graph with vertices labeled $1,2,\ldots, n$
with edges of the form $(i,N(i))$.  The non-zero contributions
in (\ref{hess}) can be efficiently computed using graph coloring techniques, as explained next.

\subsubsection{Graph coloring and Hessian vector products}

The Hessian-vector product of a function $f$ and vector $v \in \mathbb{R}^n$, is 
\begin{equation*}
    Hv = \nabla^2 f(x)v = \nabla[\nabla f(x)\cdot v ] = \nabla g(x), \quad g(x) \equiv \grad f(x) \cdot v
\end{equation*}
Computing this involves automatic differentiation to find the gradient of the scalar function $g$, rather than directly computing the Hessian of $f$.
The computational complexity in calculating a dense Hessian is comparable to computing the HVP over all $n$ coordinate vectors.   On the other hand,
for sparse matrices, we need to obtain only a small set of binary vectors $d_1,d_2,...,d_k,$ such that $Hd_1,Hd_2,...,Hd_k$ determine $H$ entirely. Finding a
minimal set of such vectors can be reformulated as a graph coloring problem, and particularly for symmetric matrices, a star coloring problem \cite{Gebremedhin2005-pq}.
In contrast to usual graph coloring, star coloring requires every path of length three to use three distinct colors. 
Vectors $d_j$ are constructed so their non-zero components are those with color $j$, and the non-zero components of the HVP $H d_j$ correspond to a subset of Hessian matrix entries.  Star coloring is an NP-hard problem,  so in practice we utilize an approximate
star coloring algorithm \cite{Coleman1984-hd,Powell1979-sx} which has nearly optimal performance.

\subsection{Numerical Comparison of Algorithms}
  \label{sec:comp}
To motivate our choice of method, we illustrate the behavior of several first- and 
second order algorithms. 
These include L-BFGS, a limited-memory quasi-Newton method \cite{Liu1989-we},
nonlinear conjugate gradients \cite{Fletcher1964-kx}, 
trust region newton conjugate gradient method (trust-ncg)\cite{Steihaug1983-xq},
trust-krylov, a trust-region method using Krylov subspace approximations,
and AdaBelief \cite{AdaBelief}, an adaptive optimizer originally developed for neural networks.

\begin{figure}
\centering
\includegraphics[height = 6cm]{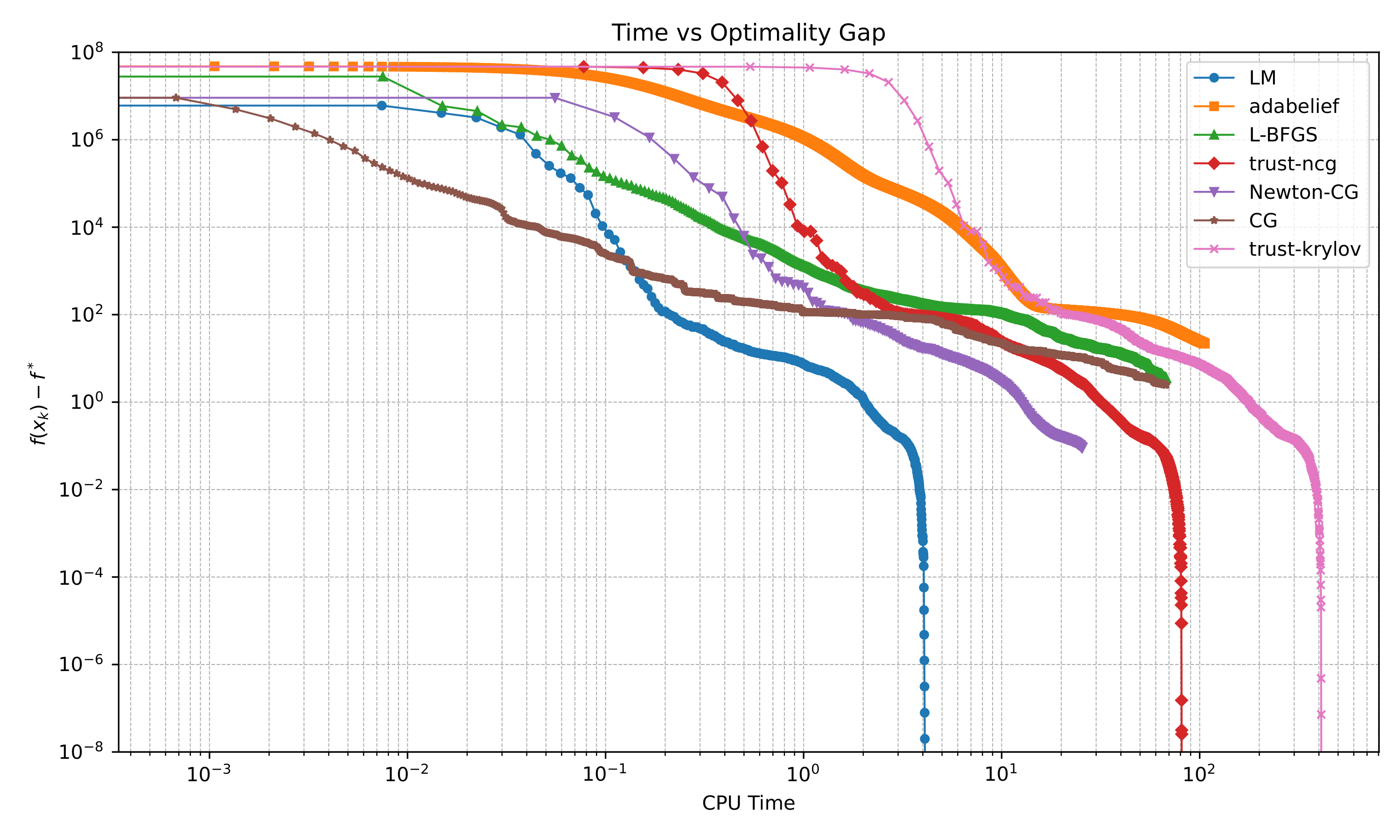}
\caption{Comparison of CPU time (seconds) between the proposed LM method and other standard optimization algorithms.}
\label{fig:OptimizationLorenz}
\end{figure}

Our tests use the objective for the problem described in section~\ref{sec:Lorenz}, 
for parameters $(\lambda, R) = (10^{-2}, 10^{-4})$.
Figure \ref{fig:OptimizationLorenz} shows CPU time versus the objective function for each method.  In this case, the generalized LM method converges at least 10 times faster that 
than the other algorithms tested.  In addition, all methods converge to the same local minimum.  While we cannot claim that the method we choose would be superior in all settings, it appears to provide reliable performance for all of the test problems 
discussed next.

\section{Results}
    \label{sec:results} 
Our system identification algorithm was tested with several well-known model differential equations. Synthetic data were generated using an adaptive Runge–Kutta $(4,5)$ method (RK45) with absolute and relative tolerances of $10^{-12}$, and sampled at uniform intervals corresponding to the sampling rates specified in each example. Various amounts of normally distributed noise were added according to 
$$
\hat{u}_i = u_i^* + \mathcal{N}(0, \sigma^2), \quad
\sigma = \text{(\% noise) $\times$ standard deviation of $u_i^*$},
$$
where $u^*$ denotes noise-free simulation data.   We also investigated the effects of
incomplete data, by either removing continuous regions (see section \ref{sec:VanderPol}), or removing data randomly  (see section \ref{sec:Lorenz}).  Note that the latter case could equivalently represent non-uniform data.  

To generate reliable initial states for optimization, we linearly interpolate missing state values, but otherwise
perform no further pre-processing. To initialize parameters we set $\theta_k = 1$ and normalize the library functions by setting $\tilde{N}_k(u) = N_k(u)/\|N_k(\hat{u})\|_2$, ensuring that each term is initially scaled to a comparable magnitude before optimization.  All results are obtained using a grid search over regularization weights $(\lambda, R) = (10^i, 10^j)$ with $i \in \{-3, -2, \ldots, 3\}$ and $j \in \{-4, -3, \ldots, 0\}$, resulting in 35 candidate models. For validation, we masked out every third time step, which was used to find the best regularization weight pair as described in section \ref{sec:Validation}.

Evaluation of system identification methods is performed using both measures of the accuracy
of fit and the success rate in identifying the model structure.
The first category of metrics uses the relative parameter and state errors
\begin{equation}
    RE(\theta) = \frac{\|\theta - \theta^*\|_{F}}{\|\theta^*\|_{F}}, \quad RE(u) = \frac{\|u - u^*\|_{F}}{\|u^*\|_{F}}.
\end{equation}
where $\|\cdot\|$ is the Frobenius norm and $\theta^*, u^*$ are the known ground truth parameters and states.  A measure of the ability to correctly identify terms is given by the 
true positivity ratio \cite{TPR_Lagergren,messenger2021weak2}
\begin{equation}
    \text{TPR}(\theta) = \frac{\text{TP}}{\text{TP + FN + FP}}.
\end{equation}
Here, TP denotes the number of active terms which are correctly identified, FN is the number of terms incorrectly set to zero, and FP is the number of terms which are incorrectly nonzero. 

Our approach is compared to the popular WSINDy method \cite{messenger2021weak2},
which performs sparse regression on a weak form of the equations.  WSINDy utilizes a predefined library of candidate equation terms
as we do here, and to date offers some of the best performance for sparse model discovery with highly noisy data.    The results presented utilize publicly available codes for WSINDy \cite{messenger2021weak2}, used here without any customizations.   

Having a sufficiently rich library of terms based on expert knowledge or known properties like symmetries would be typical in practice, but is not always available.   In our experiments, we provide candidate libraries composed of polynomial terms, up to some predefined order.
In the following examples, the true model terms are a subset of our provided libraries; we explore the consequences of insufficient libraries in Appendix \ref{Appx:InsufficLibrary}. 

All data and codes used in this manuscript are publicly available on GitHub at \url{https://github.com/TeddyMeissner/DEDi}.

\subsection{Van der Pol Oscillator} \label{sec:VanderPol}
The Van der Pol system 
\begin{align}
\label{Eqn: VanDer}
    \dot{x} &= y, \\
    \notag \dot{y} &= \mu(1-x^2)y - x, 
\end{align}
arising from a model of electrical circuits, is a non-conservative oscillator with non-linear damping. It is widely used in system identification and nonlinear dynamics studies, offering a classic example of limit cycle behavior and self-sustained oscillations.

To generate data, we used $\mu = 2$ and initial condition $(x_0, y_0) = (0, 2)$.
Simulations were performed over the interval $t \in (0, 10)$, and the resulting trajectories were sampled at uniform intervals of $\Delta \hat{t} = 0.02$,  resulting in a total of 1002 data points across 501 time steps. 

We report results using a polynomial library of degree $p = 3$ and $p = 6$.
For $p=3$, the library is
\begin{equation*}
\label{Vanderpol-Library}
N =
\begin{bmatrix}
\dot{x}_1 \\
\dot{x}_2 
\end{bmatrix} - 
\begin{bmatrix}
x_1 & x_2  & x_1 x_2 & x_1^2 & x_2^2 & x_1^2x_2 & x_1 x_2^2 & x_1^3 & x_2^3 \\
x_1 & x_2  & x_1 x_2 & x_1^2 & x_2^2 & x_1^2x_2 & x_1 x_2^2 & x_1^3 & x_2^3 
\end{bmatrix}
\bigcirc \theta,
\end{equation*}
where $\theta \in \mathbb{R}^{2 \times 9}$ contains 18 coefficients to be optimized while $p=6$ results in 54 coefficients.

\begin{figure}
    \centering
    \includegraphics[width=1.0\linewidth,trim=5 0 10 0, clip]{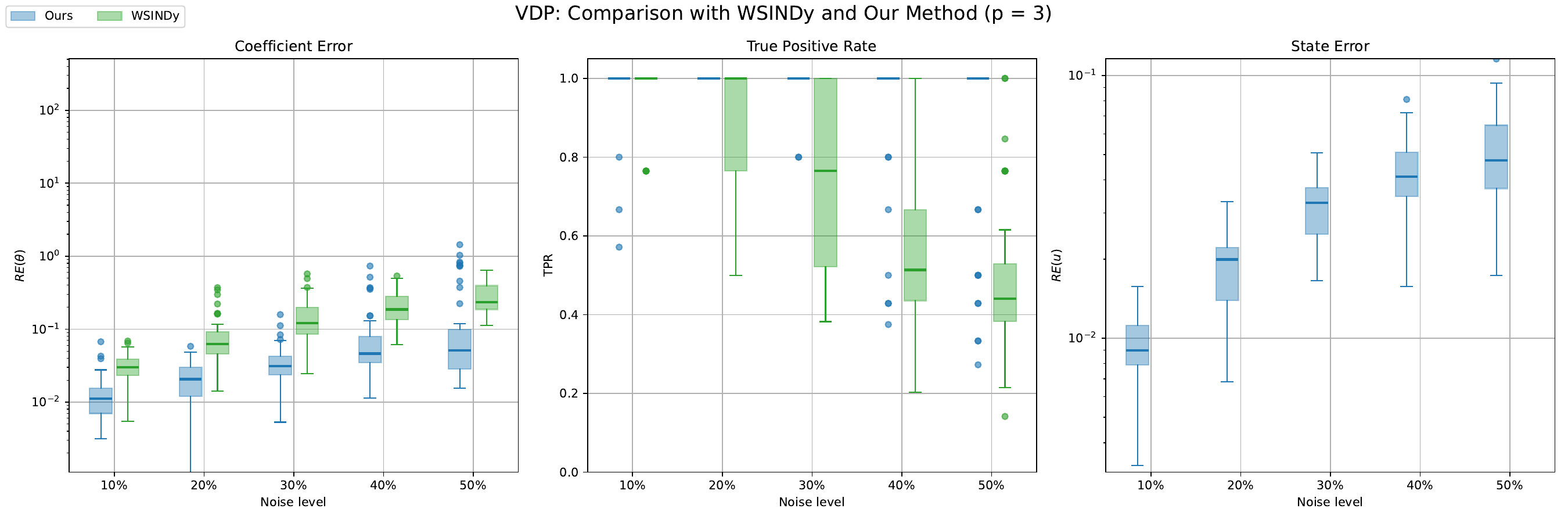} \includegraphics[width=1.0\linewidth,trim=05 0 10 0, clip]{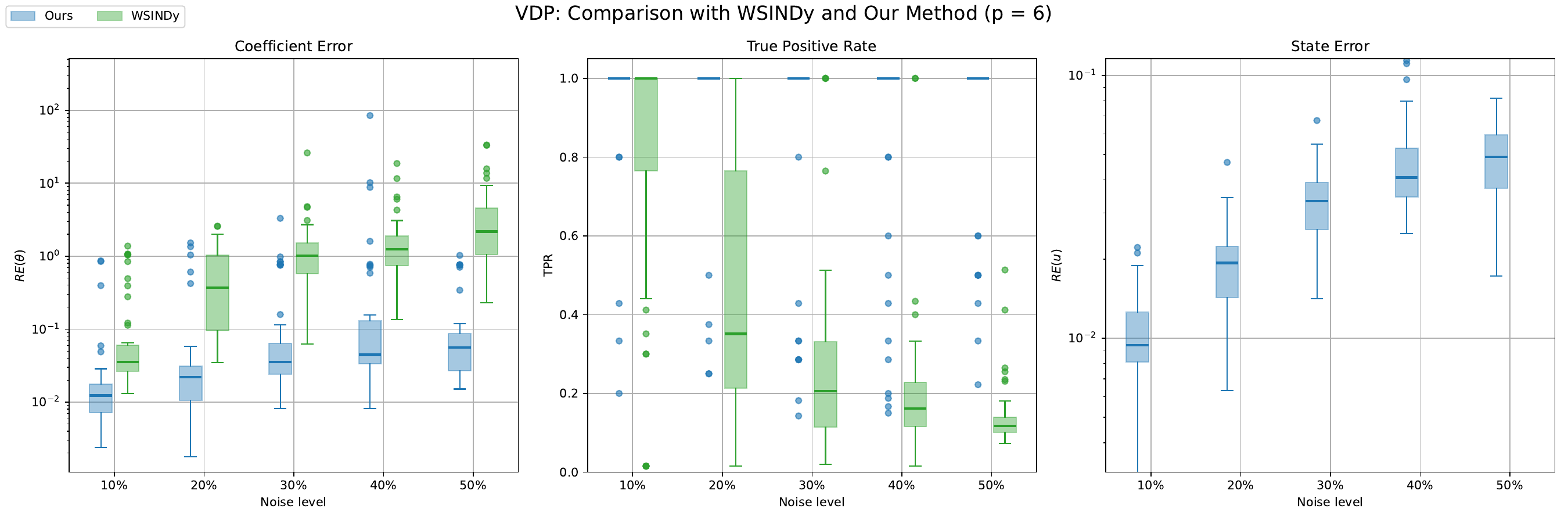}
    \caption{Comparison of our algorithm vs. WSINDy for the Van der Pol system, for monomial libraries with maximum
    orders $p=3,6$ with $18$ and $54$ terms, respectively.        
    Each boxplot summarizes 50 independent noise realizations per noise level. A line at TPR = 1 indicates that most realizations achieve perfect recovery for our method.}
\label{fig:VDP_WSINDy_Comparison_boxplot}
\end{figure}

Figure \ref{fig:VDP_WSINDy_Comparison_boxplot} shows a comparison between our approach and WSINDy using boxplots of relative parameter error $RE(\theta)$ and true positive rate (TPR) across 50 independent noise realizations at each noise level. For each boxplot, the central line denotes the median, the box spans the interquartile range (25th–75th percentile), the whiskers extend to the most extreme non-outlier values, and individual points beyond the whiskers represent outliers.

We note our approach consistently outperforms WSINDy across all noise levels, most notably at high noise. In addition,
these results demonstrate that our method remains effective when applied to a
larger candidate library. The state recovery error $RE(u)$
(which cannot be compared directly to WSINDy without integration of its discovered system), remains remarkably small despite large observation noise. 

\begin{figure}
\centering
\includegraphics[width = 15cm]{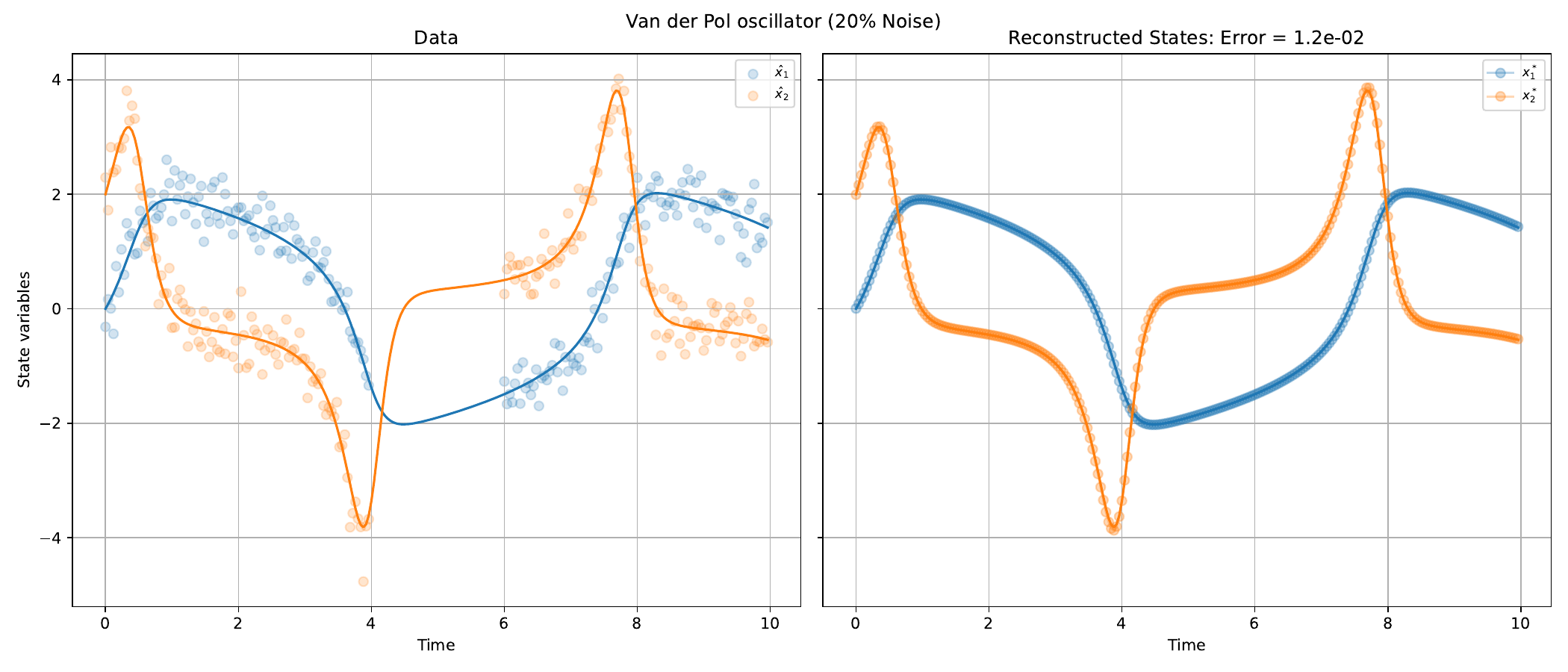}
\caption{Visualization of the Van der Pol oscillator with a missing continuous region of data. Left: the true state of the system (solid), and observed data (dots);  Right:
the reconstructed state variables.}
\label{fig:VanDerNumerical}
\end{figure}

\begin{figure}
    \centering
    \includegraphics[width=1.0\linewidth,trim=5 0 10 0, clip]{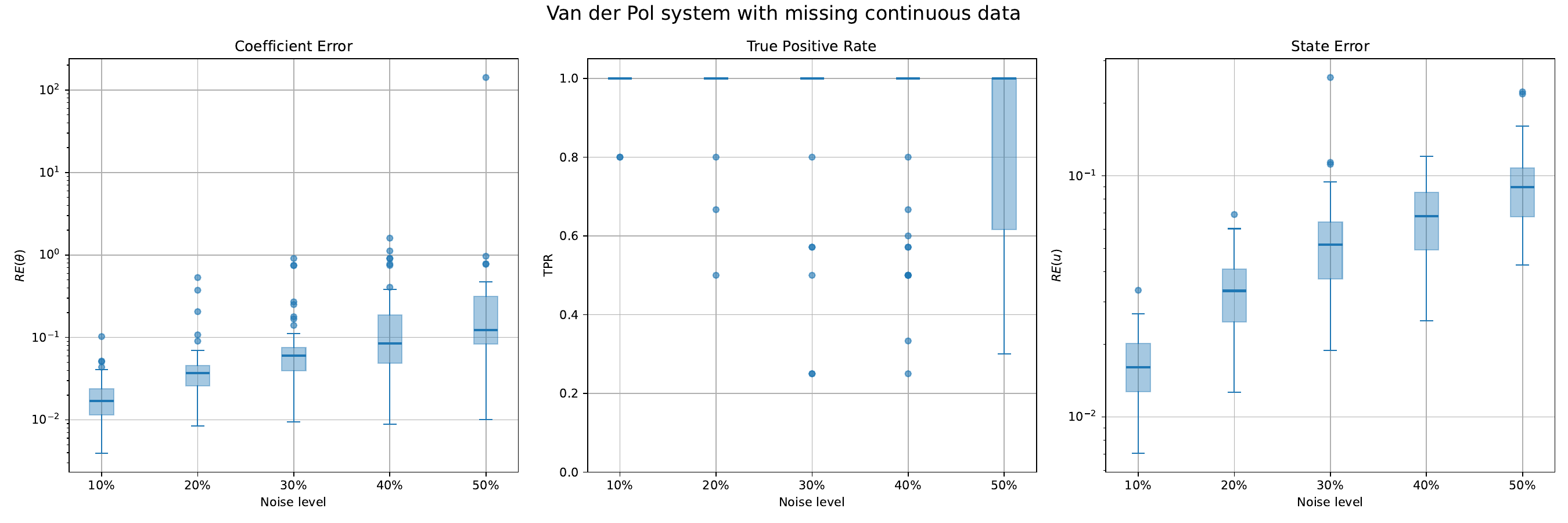}
    \caption{Each boxplot summarizes 50 independent noise realizations per noise level for the VDP system with missing continuous data.}
    \label{fig:VDP_cont_missing_boxplots}
\end{figure}

A more extreme scenario with both measurement noise and continuous gaps in the observations was considered using a polynomial library of degree $p=3$. (see figure \ref{fig:VanDerNumerical}). The data shown here sampled at $\Delta \hat{t} = .04$, and data between $t\in (4,6)$ are completely removed. This results in a total of 402 data points across 201 time steps. Figure \ref{fig:VDP_cont_missing_boxplots} shows boxplots of relative parameter error and true positive rate across 50 samples at varying levels of noise. These results emphasize the robustness of our method under challenging conditions.

Although the system's structure can usually be determined, occasionally spurious terms appear, mostly at high noise levels.  At $10\%$ noise, only two in fifty trials 
misidentifies model terms, and a typical failure mode is provided by the following example:
\begin{align*}
    \dot{x} &= \color{blue}{1.02 y} \color{red}{ \ + \ .03 x}, \\
    \dot{y} &= \color{blue}{-1.13 x + 2.00 y - 2.09 x^2 y},
\end{align*}
which contains an extra, but fairly negligible, term. At $50\%$ noise, 33 out of 50 results were
structurally correct, although more severe misidentification problems are occasionally seen.
For instance, in one situation the algorithm prematurely discarded a key term, leading to
\begin{align*}
    \dot{x} &= \color{blue}{1.05 y}, \\
    \dot{y} &= \color{blue}{2.13 y - 1.48 x^2 y} \color{red}{\ - \ 0.19 y^3 - 0.25 x y^2}.
\end{align*}

Despite this, the method consistently captures the core dynamics across most trials, even in the presence of substantial data loss and noise, highlighting its robustness and reliability.

\subsection{The Lorenz system} \label{sec:Lorenz}
The Lorenz equations \cite{Lorenz1963-tp}
\begin{align}
\label{Eqn: Lorenz}
    \dot{x} &= \sigma(y-x) \\
    \notag \dot{y} &= x(\rho-z) + y \\
    \notag  \dot{z} &= xy - \beta y
\end{align}
were originally derived from models of atmospheric convection,
and have become the standard example of chaotic behavior in ordinary differential equations.   They have also emerged as a canonical test case in studies of system identification due to their deterministic yet complex behavior.

Simulations used $(\sigma, \rho, \beta) = (10, 28, 8/3)$ and initial conditions $(x_0, y_0, z_0) = (-8, 8, 27)$.
The system was simulated over the interval $t \in (0, 10)$, and trajectories were sampled at uniform intervals of $\Delta \hat{t} = 0.02$, resulting in a total of 1503 data points across 501 time steps. The library of equation terms consisted of monomials up to degree $p=3$ and $p=4$. For $p=3$, the library is

\begin{equation*}
\label{Lorenz-Library}
N =
\begin{bmatrix}
\dot{x}_1 \\
\dot{x}_2 \\
\dot{x}_3
\end{bmatrix} -
\begin{bmatrix}
x_1 & \cdots & x_1 x_2 & \cdots & x_3^2 & x_1 x_2 x_3 & \cdots & x_3^3 \\
x_1 & \cdots & x_1 x_2 & \cdots & x_3^2 & x_1 x_2 x_3 & \cdots & x_3^3 \\
x_1 & \cdots & x_1 x_2 & \cdots & x_3^2 & x_1 x_2 x_3 & \cdots & x_3^3 \\
\end{bmatrix}
\bigcirc \theta,
\end{equation*}
where $\theta \in \mathbb{R}^{3 \times 19}$ is a coefficient matrix containing 57 entries. For $p = 4$, the resulting library contains 102 terms.

\begin{figure}
    \centering
    \includegraphics[width=01.0\linewidth,trim=5 0 10 0, clip]{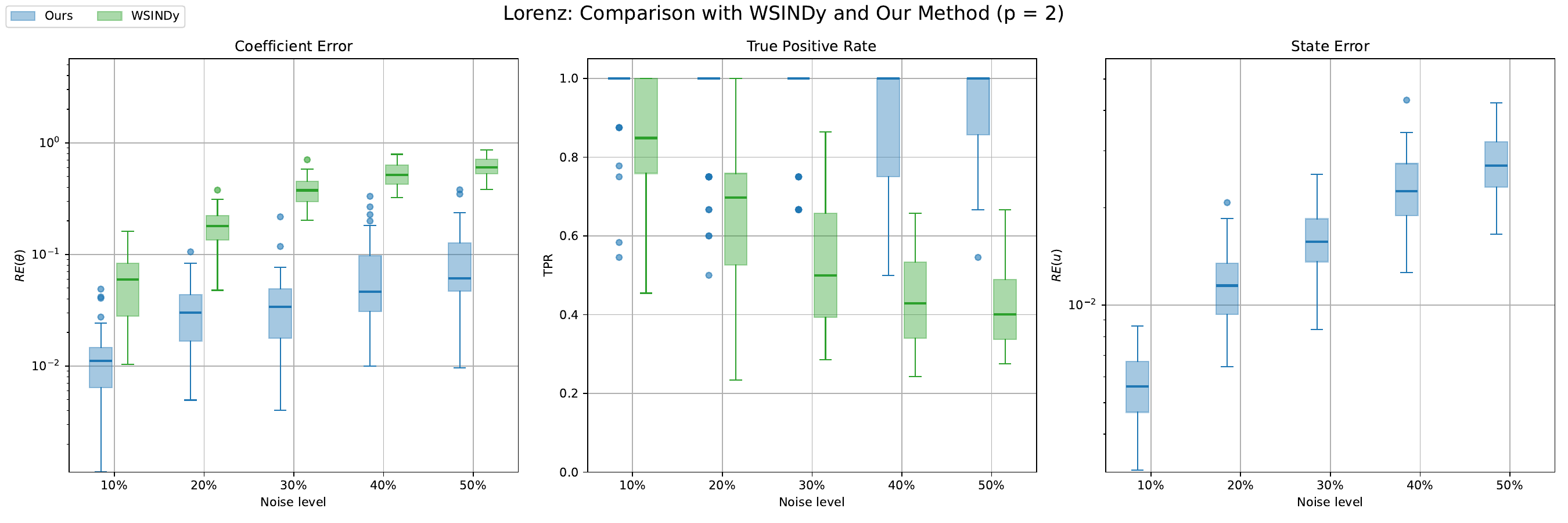} \\
    \includegraphics[width=01.0\linewidth,trim=5 0 10 0, clip]{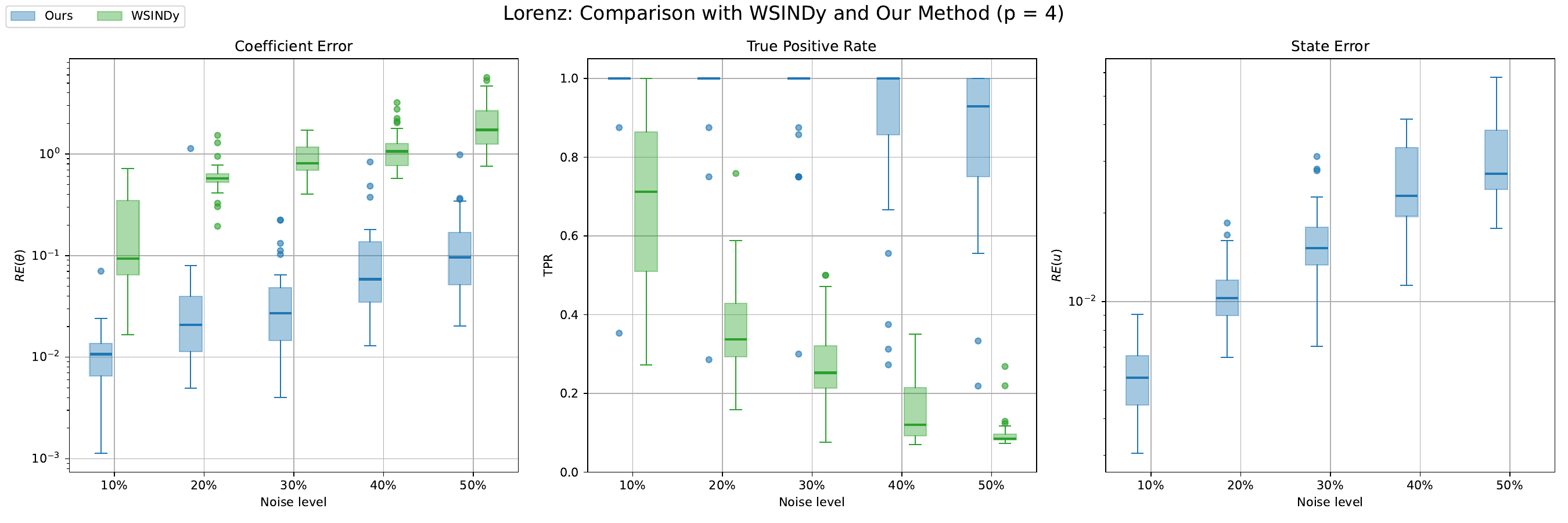}
    \caption{Comparison of our algorithm vs. WSINDy for the Lorenz system, for monomial libraries with maximum
    degrees $p=3,4$ with $57$ and $102$ terms, respectively. Each boxplot summarizes 50 independent noise realizations per noise level.}
    \label{fig:Lorenz_WSINDy_Comparison_boxplots}
\end{figure}

In figure \ref{fig:Lorenz_WSINDy_Comparison_boxplots}, we again benchmark our approach using boxplots of relative parameter errors and recovery rates across 50 independent realizations for each noise level.  The results are similar to the previous example,
with both better parameter estimates and recovery rates compared to WSINDy.
Even with 50\% noise, it was often possible to recover an accurate representation of
the ground truth model.

\begin{figure}
\centering
\includegraphics[width = 15cm]{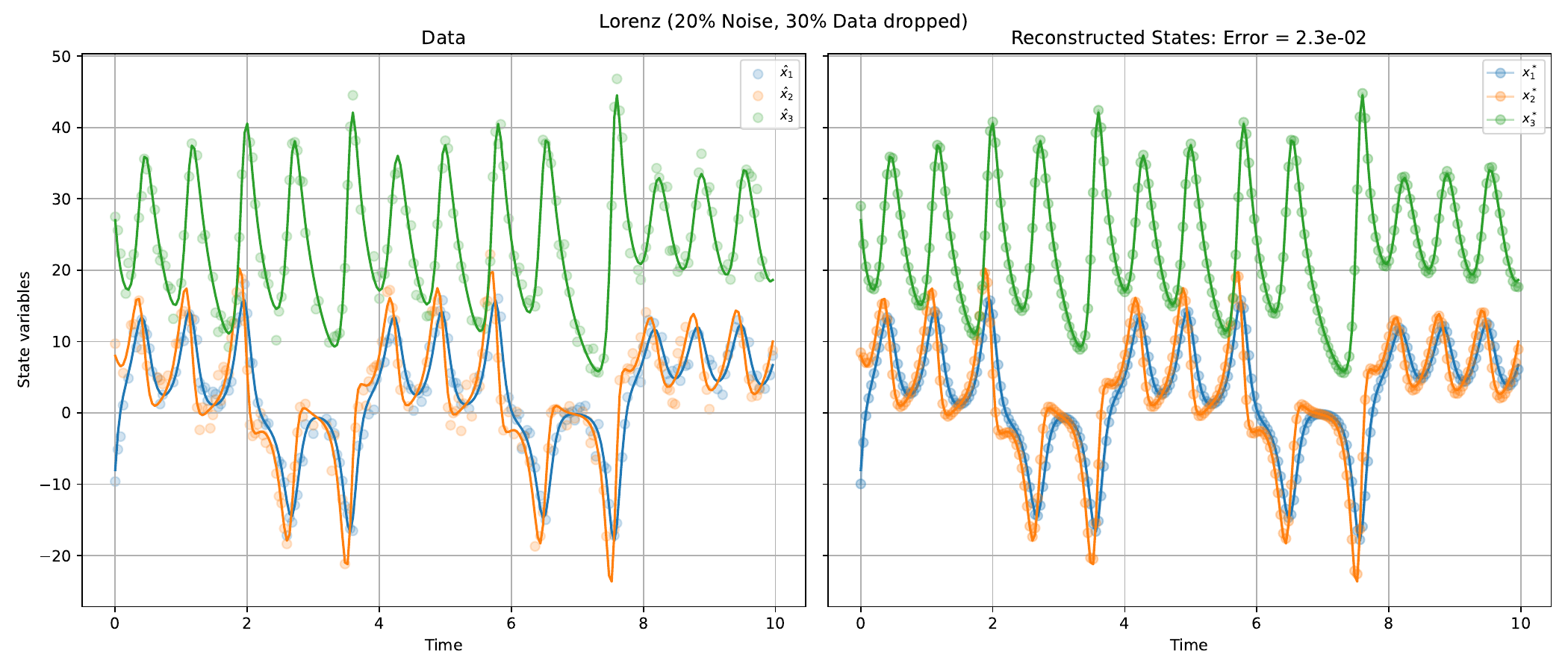}
\caption{Visualization of the Lorenz system with 20\% noise and 30 \% randomly dropped data. 
Left: the true state of the system (solid), and observed data (dots);  Right:
the reconstructed state variables.}
\label{fig:Lorenz_rand_missing}
\end{figure}

\begin{figure}
    \centering
    \includegraphics[width=1.0\linewidth,trim=5 0 10 0, clip]{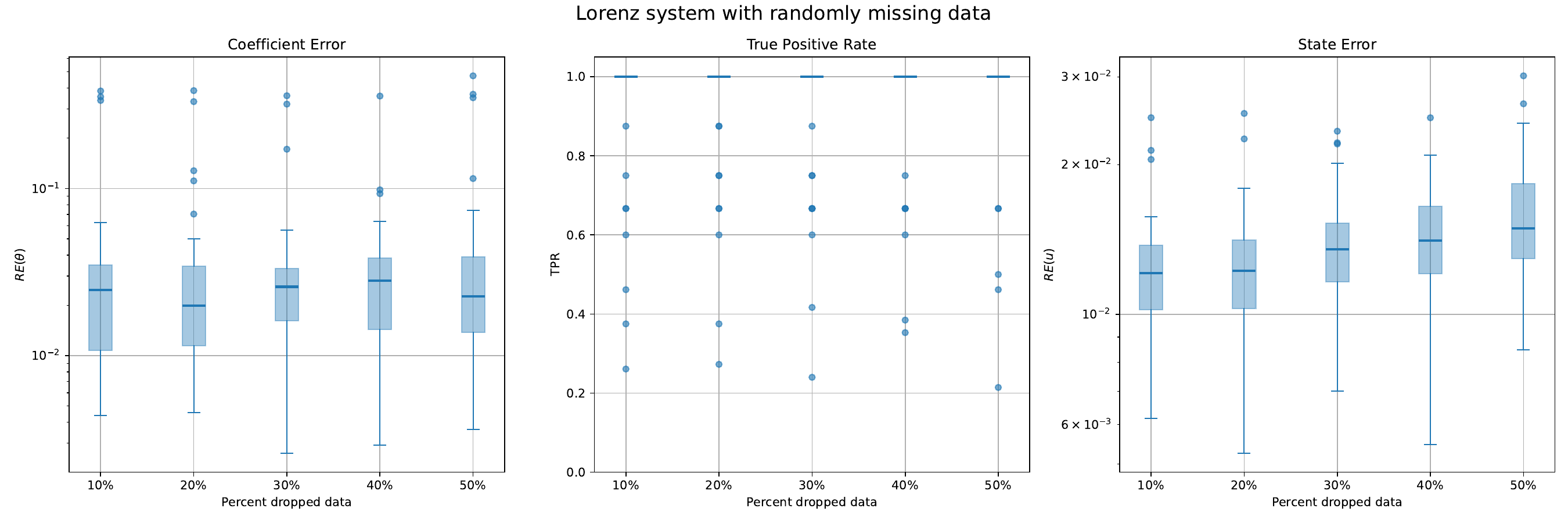}
    \caption{Each boxplot summarizes 50 independent noise realizations per noise level for the Lorenz system with various levels of randomly dropped data.}
    \label{fig:Lorenz_random_drop_boxplot}
\end{figure}

To assess performance under incomplete observations, we added 20\% noise and randomly removed data points from each state component. Results were performed with a polynomial library of degree $p = 3$. 
Figure \ref{fig:Lorenz_rand_missing} presents a representative case with 30\% missing data, demonstrating strong agreement between the recovered and true state trajectories. Figure \ref{fig:Lorenz_random_drop_boxplot} shows the
resulting parameter and state errors, along with the recovery rates.  
For even large quantities of missing data, model parameters can be estimated with good accuracy.

\subsection{The Lorenz-96 system} \label{sec:Lorenz96}
To evaluate our method for higher dimensional systems, we utilized the Lorenz-96 model \cite{Lorenz96}, which also exhibits chaotic dynamics, and is often invoked to evaluate
data-driven methods for dynamical systems.

For model dimension $N = 5$, the dynamics are governed by
\begin{equation}
\label{Eqn: L96}
\dot{x}_i = (x_{i+1} - x_{i-2})x_{i-1} - x_i + F, \quad i = 1, \dots, 5. 
\end{equation}
A standard forcing value of $F = 8 $ is used in our simulations with an initial condition slightly perturbed from the equilibrium state $x(0) = (F+0.01,F,F,F,F)$.  Simulations were performed over the interval $t \in (0,10)$ with sampling rate $\Delta \hat{t} = 0.04$, resulting in a total of 1255 data points across 251 time steps.

The candidate library consisted of all polynomial terms up to degree 2, so that
\begin{equation}
\label{Eqn: L96-Library}
N =
\begin{bmatrix}
\dot{x}_1 \\
\vdots \\
\dot{x}_5
\end{bmatrix}
-
\begin{bmatrix}
1 & x_1 & \cdots & x_5 &
x_1^2 & x_1 x_2 & \cdots & x_5^2 \\
\vdots & \vdots & \vdots & \vdots &
\vdots & \vdots & \vdots & \vdots \\
1 & x_1 & \cdots & x_5 &
x_1^2 & x_1 x_2 & \cdots & x_5^2 \\
\end{bmatrix}
\bigcirc \theta,
\end{equation}
where $\theta \in \mathbb{R}^{5 \times 21}$ is a coefficient matrix containing 105 entries.

\begin{figure}[H]
    \centering
    \includegraphics[width=1.0\linewidth,trim=5 0 10 0, clip]{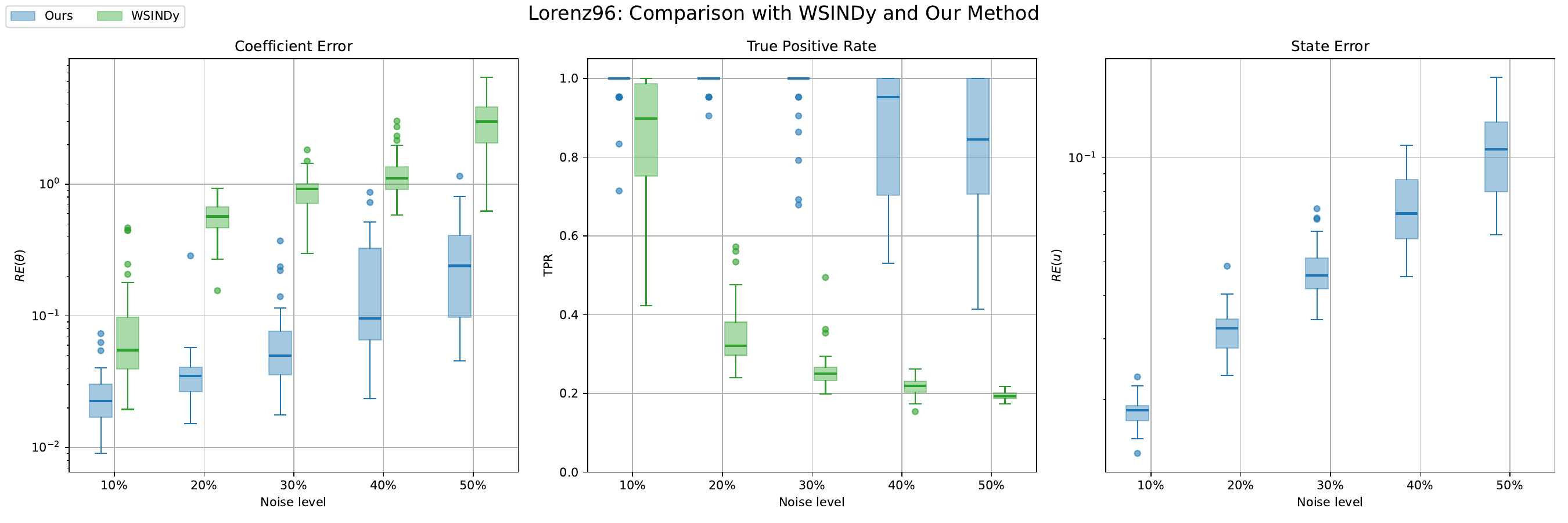}
    \caption{Comparison of our algorithm vs. WSINDy. Each boxplot summarizes 50 independent noise realizations per noise level.}
    \label{fig:Lorenz96_WSINDy_Comparison_boxplots}
\end{figure}

\begin{figure}[H]
\centering
\includegraphics[height = 6cm]{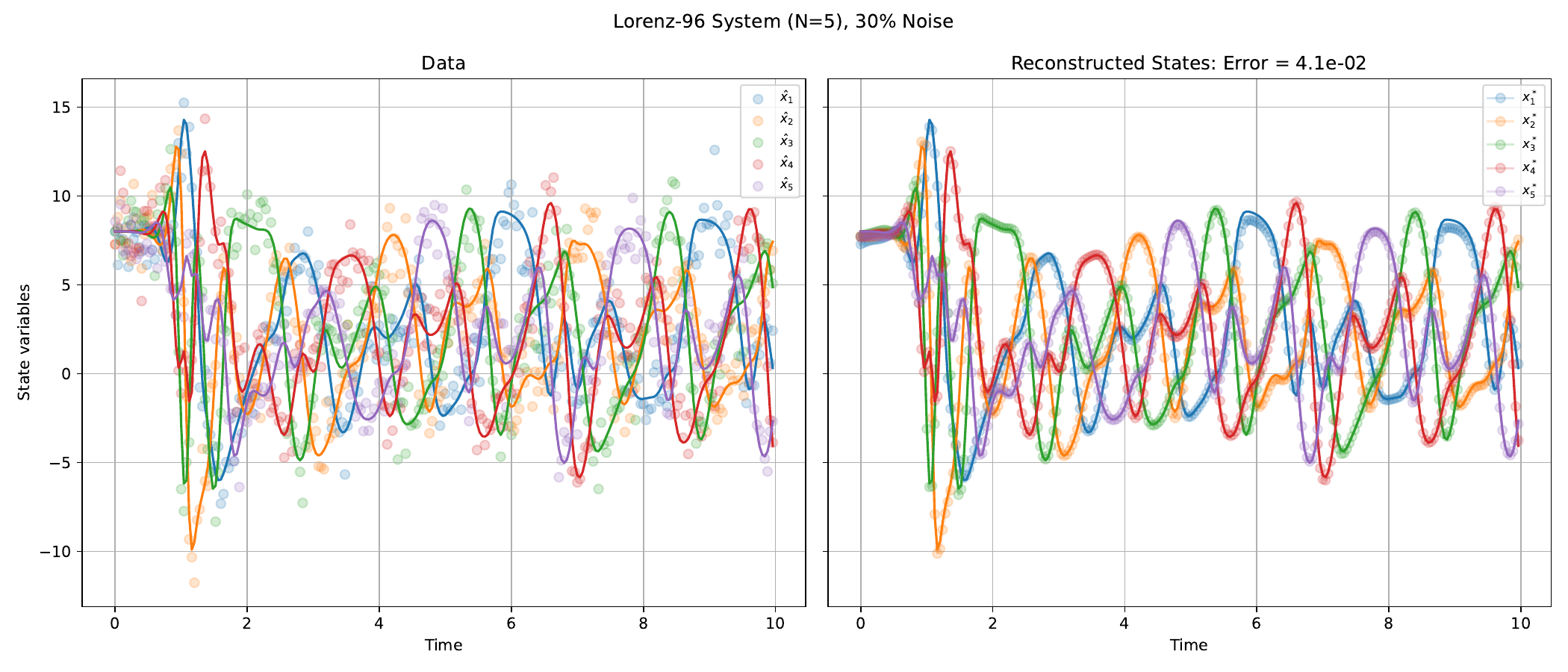}
\caption{Visualization of the Lorenz96 system's components over time. 
Left: the true state of the system (solid), and observed data (dots);  Right:
the reconstructed state variables. }
\label{fig:Lorenz96}
\end{figure}

We compare our method to WSINDy in figure \ref{fig:Lorenz96_WSINDy_Comparison_boxplots}, where we see again the ability of our method to have significantly higher model recovery rates, especially as the level of noise increases.  Despite the large candidate library, the median TPR stays above 0.9 for noise levels up to 40\%.  Figure \ref{fig:Lorenz96} depicts an example of state reconstruction under 30\% noise, along with the corresponding discovered model in \eqref{eq:Lorenz96_reconstruct_equation}.

\begin{equation}\label{eq:Lorenz96_reconstruct_equation}
\begin{aligned}
\dot{x}_1 &= 
8.005 
- 0.983 x_1 
+ 1.010 x_2 x_5 
- 1.016 x_4 x_5 \\
\dot{x}_2 &= 
8.144 
- 1.032 x_2 
+ 1.021 x_1 x_3 
- 1.035 x_1 x_5 \\
\dot{x}_3 &= 
7.840 
- 0.954 x_3 
- 0.999 x_1 x_2 
+ 1.017 x_2 x_4 \\
\dot{x}_4 &= 
8.020 
- 1.042 x_4 
- 1.015 x_2 x_3 
+ 1.045 x_3 x_5 \\
\dot{x}_5 &= 
7.763 
- 0.994 x_5 
+ 0.997 x_1 x_4 
- 1.009 x_3 x_4 \\
\end{aligned}
\end{equation}

\subsection{Nonlinear term libraries and the Colpitts oscillator} \label{Sec:Colpitts}

One of the current method's virtues is that candidate models are not required to be
linear in parameters, unlike other sparse regression methods that rely on linear least squares.  To illustrate this feature, we study the
Colpitts oscillator as described in \cite{Colpitts_Abarbanel2013} and given by,
\begin{equation}
\label{Eqn: Colpitts}
\begin{aligned}
    \dot{x} &= \alpha_C z \\
     \dot{y} &= \eta (1 - \exp(ax) + z)\\
     \dot{z} &= -\gamma(x+y) - qz. 
\end{aligned}
\end{equation}
where the parameter $a$ in the term $\exp(ax)$ is also unknown {\it apriori}.
Parameters in this class are estimated in the regression method, although are ignored for purposes of sparsity (the prefactor of such a term, on the other hand, could be set to zero and
the term removed if the algorithm permits).  Ribera et al.\ \cite{ribera2022model} similarly observe
that such library terms may be included, and demonstrated successful recovery of this system 
under noise-free conditions.  Our aim here is to show even in the presence of large noise, the correct system can be recovered.

Simulated data used $(\alpha_C, \eta, \gamma, q, a) = (5, 6.2723, 0.0797, 0.6898, -1)$ and initial condition $(x_0,y_0, z_0) = (0.01,0,0)$, and were generated for the interval $t \in (0,50)$ with sampling rate $\Delta \hat{t} = 0.1$, resulting in a total of 1500 data points across 500 time steps. 
The initial library provided was: 
\begin{equation*}
\label{Colpitts-Library}
N =
\begin{bmatrix}
\dot{x}_1 \\
\dot{x}_2 \\
\dot{x}_3
\end{bmatrix}
-
\begin{bmatrix}
1 & x_1 & x_2 & x_3 & x_1 x_2 & x_1 x_3 & x_2 x_3 & x_1^2 & x_2^2 & x_3^2 \ \ e^{a x_1}  \\
1 & x_1 & x_2 & x_3 & x_1 x_2 & x_1 x_3 & x_2 x_3 & x_1^2 & x_2^2 & x_3^2 \ \ e^{a x_1}  \\
1 & x_1 & x_2 & x_3 & x_1 x_2 & x_1 x_3 & x_2 x_3 & x_1^2 & x_2^2 & x_3^2 \ \ e^{a x_1}  \\\end{bmatrix}
\bigcirc \theta,
\end{equation*}
where $\theta \in \mathbb{R}^{3 \times 12}$ is a linear coefficient matrix, thus there are $37$ coefficients to optimize over including the nonlinear coefficient $a$. Due to the large difference in scales between system components, we normalize the states as $\bar{x}_i = \hat{x}_i / \text{std}(\hat{x}_i)$ and perform optimization over $\bar{x}$. The original system coefficients are then recovered by appropriately rescaling the optimized coefficients. 

\begin{figure}
    \centering
    \includegraphics[width=1.0\linewidth,trim=5 0 10 0, clip]{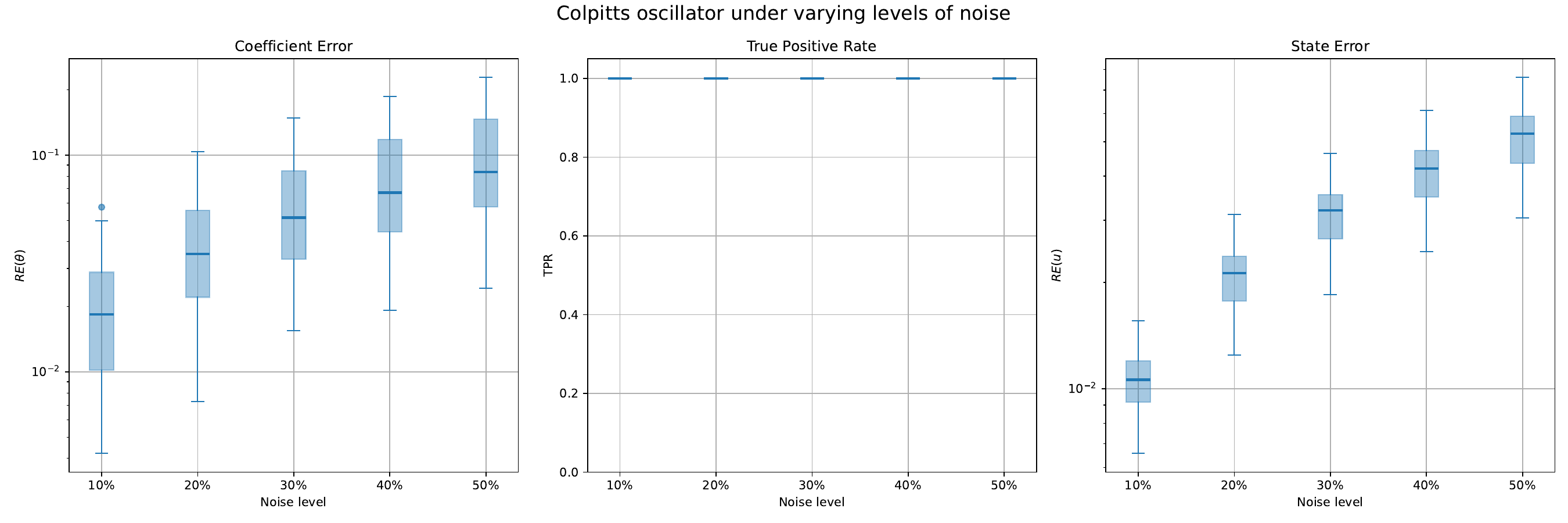}
    \caption{Each boxplot summarizes 50 independent noise realizations per noise level. A line at TPR = 1 indicates most realizations achieve perfect recovery.}
    \label{fig:Colpitts_boxplots}
\end{figure}

\begin{figure}
\centering
\includegraphics[width = 15cm]{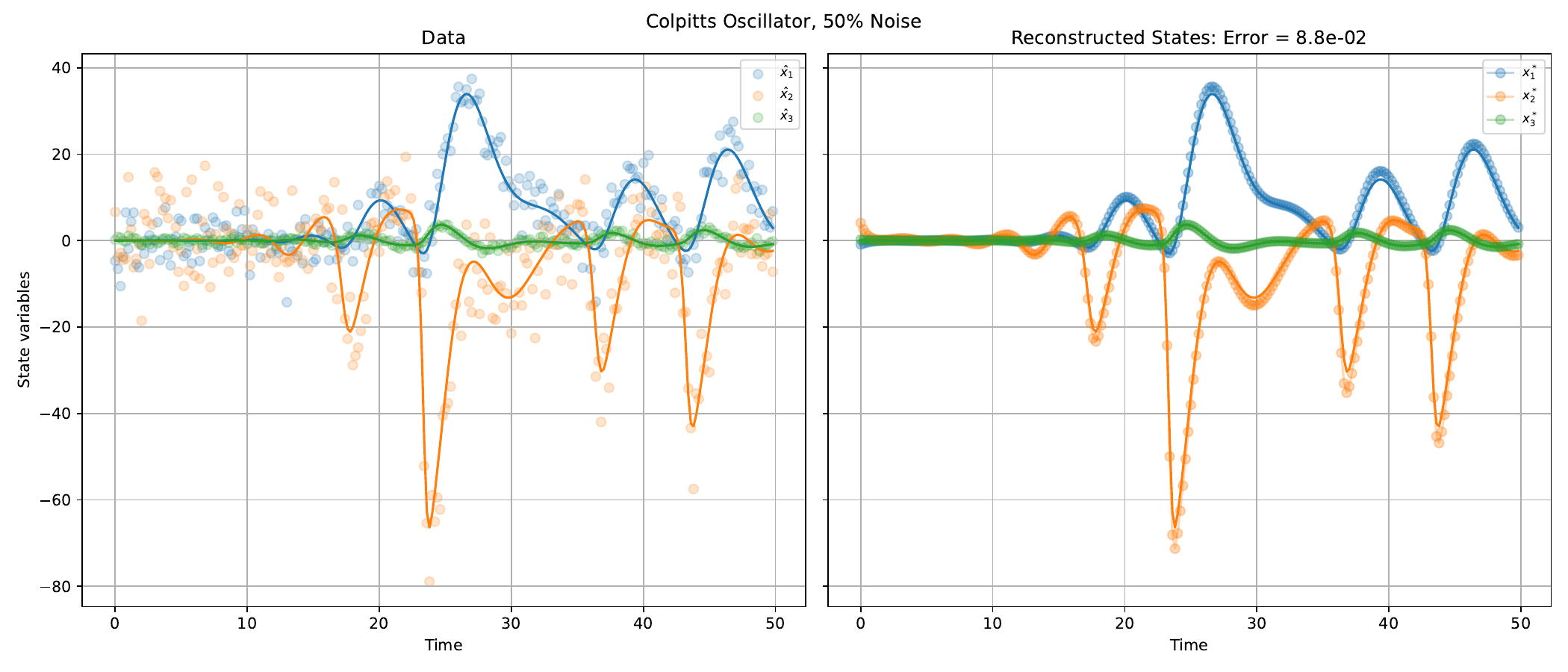}
\caption{
Visualization of the Colpitts oscillator's components over time. 
Left: the true state of the system (solid), and observed data (dots);  Right:
the reconstructed state variables.  }
\label{fig:ColpittsNoise}
\end{figure}

Figure~\ref{fig:ColpittsNoise} shows a representative example with 50\% noise added with the corresponding discovered equation 
\begin{equation}
\label{Eqn: Colpitts_discovered}
\begin{aligned}
    \dot{x} &= 5.08 z \\
     \dot{y} &= 6.68 - 5.75 \exp(-1.01x) + 6.35 z\\
     \dot{z} &= -0.08 x - 0.08 y - 0.77 z. 
\end{aligned}
\end{equation}
Figure~\ref{fig:Colpitts_boxplots} summarizes relative errors and true positive rates (TPR) from 50 trials across various noise levels. Notably, our method shows remarkable ability to consistently capture the nonlinear parameter even with substantial noise. These results underscore the potential of our approach when linear libraries are not sufficient in explaining the dynamics. 

\subsection{Model Refinement}\label{subsection:Model_Refinement}

While the previous examples focused on recovering governing equations under a fixed discretization, they assume that the data sampling rate was sufficiently fine to capture the essential system dynamics.  
Here, we relax that assumption and explore how the choice of model-grid spacing interacts with the sampling rate of the data.  
This experiment serves as a preliminary investigation toward understanding the tradeoff between discretization accuracy and model recovery in our framework.

We performed a study using the Van der Pol and Lorenz systems under 20\% measurement noise.  
Observed data were sampled at fixed intervals of $\Delta \hat{t} = 0.05, 0.1, 0.2,$ and $0.4$, while the model discretization used for optimization was independently varied between $\Delta t \in (0.025, 0.4)$.  
For each combination of sampling frequency and model grid, the discovery procedure was repeated over 20 independent noise realizations.

\begin{figure}
    \centering
    \includegraphics[width=1.0\linewidth,trim=5 0 10 0, clip]{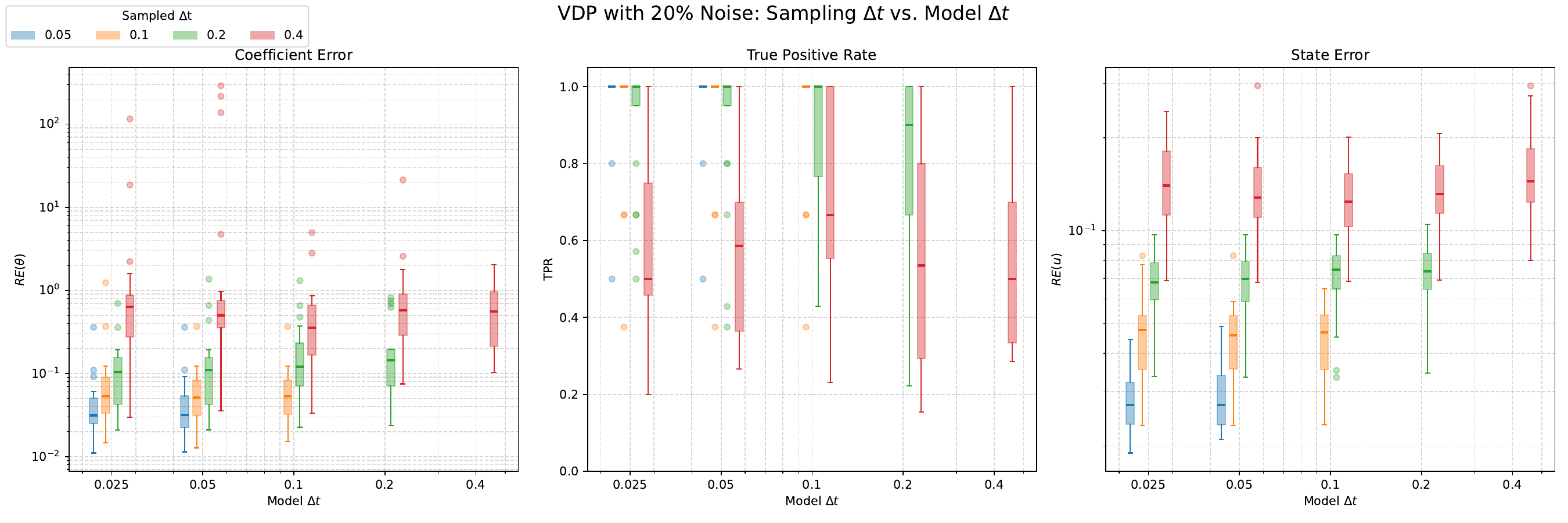}
    \includegraphics[width=1.0\linewidth,trim=5 0 10 0, clip]{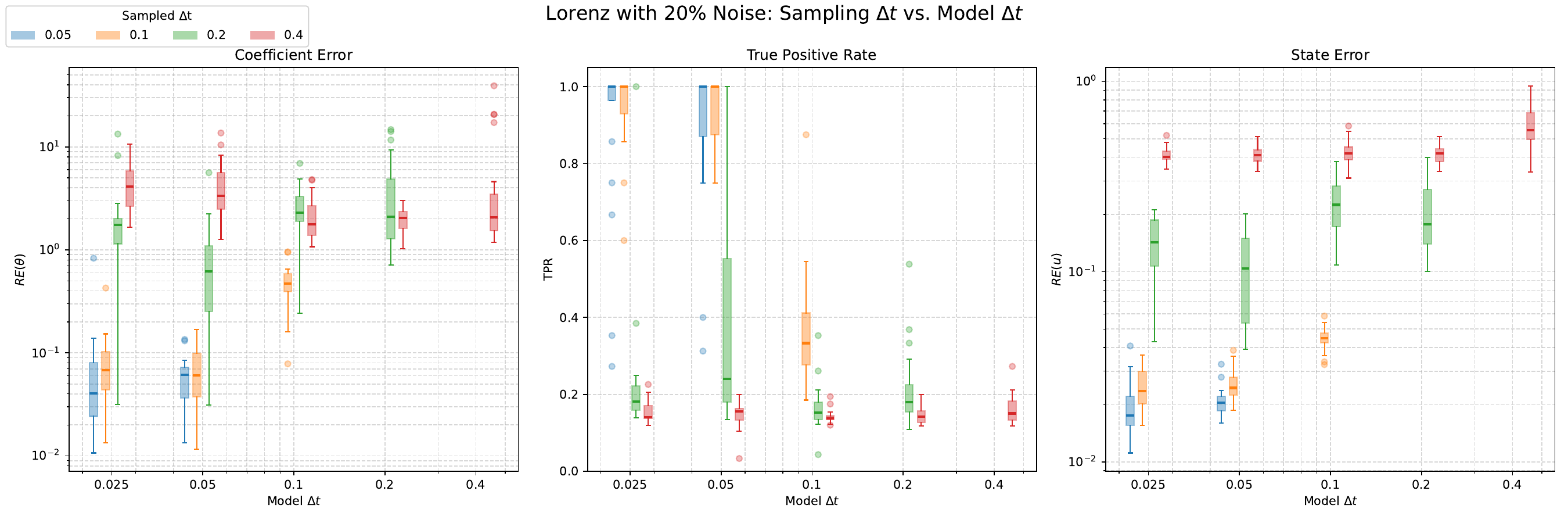}
    \caption{
    Effect of model-grid refinement on the Van der Pol system (top) and Lorenz system (bottom) under 20\% measurement noise.
    Boxplots show parameter error ($RE(\theta)$), true positive rate (TPR), and state error ($RE(u)$) as functions of model-grid spacing for different sampling intervals.}
    \label{fig:Refinement_Boxplots}
\end{figure}

Figure \ref{fig:Refinement_Boxplots} shows refinement of the model grid generally improves parameter and state accuracy up to a limit.
For the Van der Pol system, moderate refinement yields small but consistent gains, though at the coarsest sampling interval ($\Delta \hat{t} = 0.4$) refinement provides little benefit.
For the Lorenz system, refinement leads to significant improvements at finer sampling intervals ($\Delta \hat{t} \le 0.1$), but when the data are too coarsely sampled ($\Delta \hat{t} \ge 0.2$–$0.4$), the benefits saturate.
Across both systems, excessive refinement occasionally worsens performance due to the introduction of additional optimization variables.
The observed interplay provides insights into the balance between discretization bias, conditioning, and information content. We believe further work in establishing model discretization in a rigorous fashion would be a great extension to this work.

\section{Discussion and conclusions}
   \label{sec:conclusions}
   
In this study, we introduced a framework for model discovery that is robust in extreme situations of limited data and substantial noise.  With implementation of automatic differentiation,
the procedure outlined here can be extended to a wide variety of problems without significant amounts of coding.  In addition, local discretization schemes can be implemented that
allow for arbitrary accuracy of solution representations, yet provide sparsity
of the Hessian required for efficient use of second-order optimization techniques.

Our formulation was chosen to limit the requirement for significant user configuration and manual fine-tuning.  This is often a disadvantage of many data-driven algorithms where there is a need to supply numerous hyperparameters, neural architectures, as well as preprocess or filter data. 
We expect that further improvements could directly integrate hyperparameter selection 
within the main optimization algorithm.

There are a variety of natural extensions to this work.   The abstract setting given in
section \ref{sec:setup} makes no restriction to ordinary differential equations.
We expect that our framework can be adapted to applications such as
differential-algebraic systems and partial differential equations, although the latter would entail very high dimensional optimization.
We also imagine that merging our ideas with weak-form coarse discretizations
might be valuable in situations where fine-scale discretizations lead to computational intractability.  Lastly, further improvements in the underlying formulation could allow for hyperparameters to be treated as statistical parameters to be inferred, 
as in \cite{glasner2023data}.

The sophistication of our methodology is accompanied by its own set of challenges, including a complex optimization landscape and poor conditioning.   While the hybrid loss function ameliorates some
of the difficulties of non-convexity, there are nonetheless circumstances where undesirable local
minima are encountered.    Further improvement of the optimization algorithms required here
should be pursued, and will most likely require novel or bespoke methods. 

\section*{Acknowledgments}
The authors were supported through NSF awards DMS-1908968 and DMS-2407033.
Additionally, TM was supported by the Department of Defense (DoD) through the National Defense Science \& Engineering Graduate (NDSEG) Fellowship Program.

\bibliographystyle{siamplain}
\bibliography{references}

\appendix

\section{Insufficient Candidate Libraries} \label{Appx:InsufficLibrary}
A sufficiently expressive candidate library is essential for accurately recovering the underlying dynamics of a system. Here, we illustrate the result of using an insufficient library by applying our method to the Van der Pol oscillator with a basis limited to polynomials of degree at most two, thereby omitting its characteristic cubic term (see \eqref{Vanderpol-Library-insufficient}).

The data generation setup follows that of section \ref{sec:VanderPol}, with a time step of $dt = 0.04$, resulting in 201 time steps and a total of 402 data points. For model selection, we perform a grid search over the regularization parameters $(\lambda, R) = (10^i, 10^j)$ with $i \in \{-3, -2, \ldots, 3\}$ and $j \in \{-4, -3, \ldots, 0\}$, resulting in 35 candidate models, of which the winning model is shown in table \ref{tab:vdp_discovered_models}. 

\begin{figure}[H]
    \centering
    \includegraphics[width=0.42\linewidth]{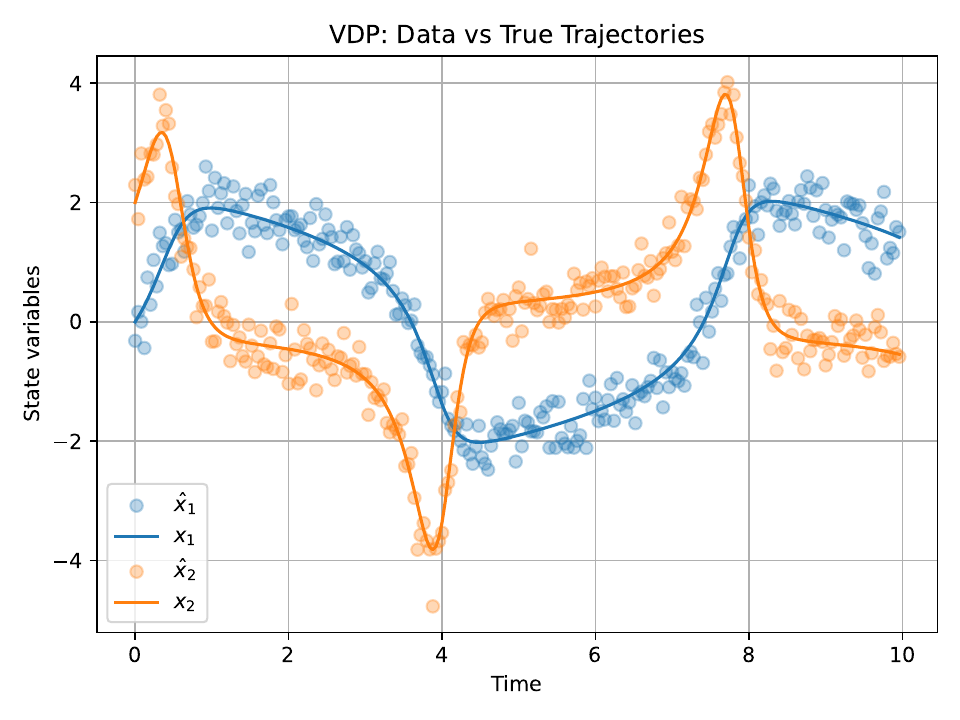}
    \includegraphics[width=0.42\linewidth]{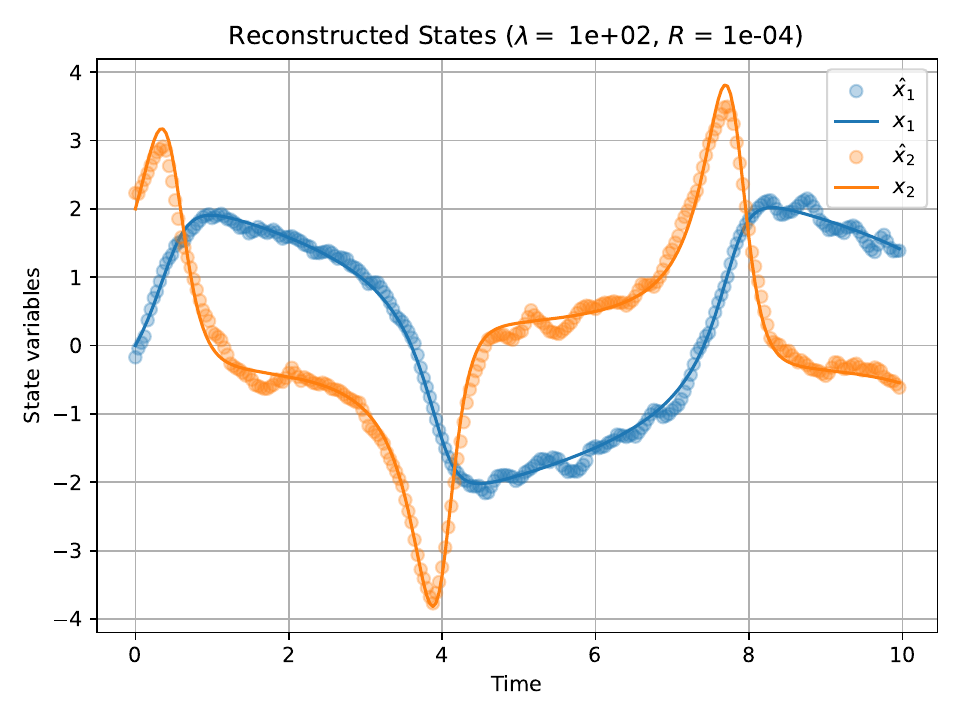}\\
    \includegraphics[width=0.42\linewidth]{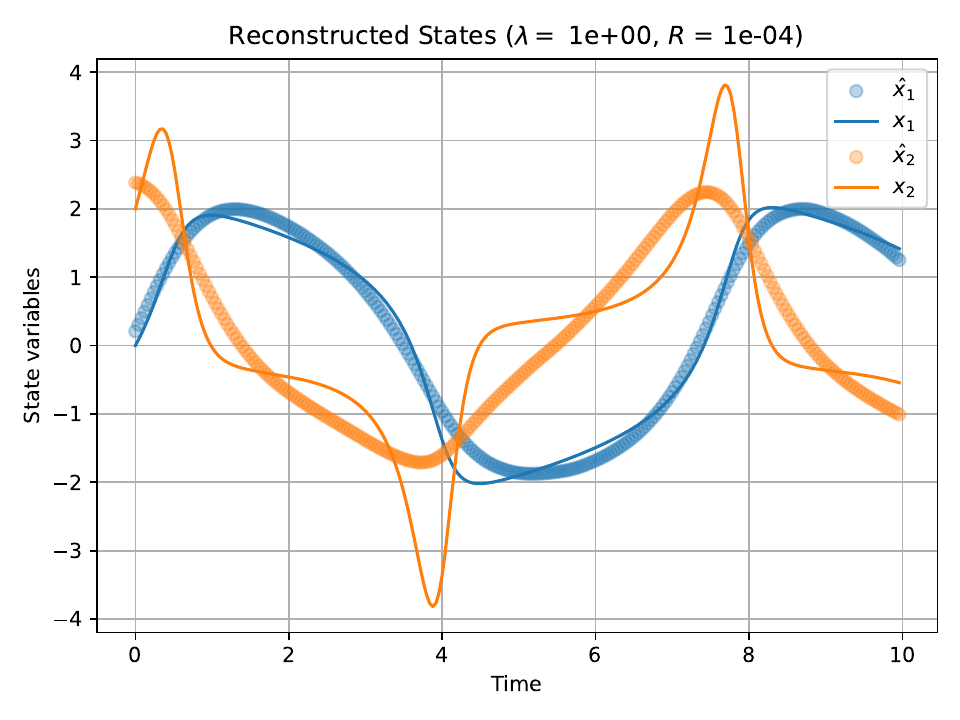}
    \includegraphics[width=0.42\linewidth]{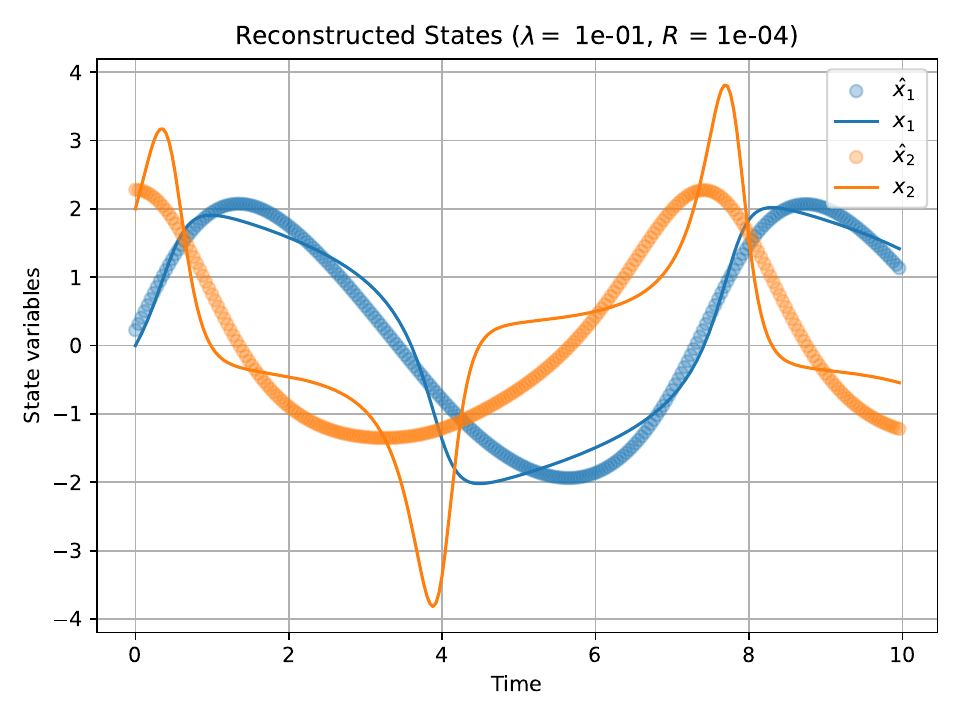}
    \caption{State estimation using recovered models for the Van der Pol oscillator using an insufficient library of polynomial terms (up to degree 2).}    \label{fig:vdp_insufficient}
\end{figure}

Among these, only three distinct model structures were present. For each structure, we retained the model with the lowest validation error. The recovered models are summarized in table \ref{tab:vdp_discovered_models} and visualized in figure \ref{fig:vdp_insufficient}. As expected, we observe a clear relationship between the regularization parameter $\lambda$ and model error: larger values of $\lambda$ correspond to higher model errors. We see here that our method selects that which is close to the data, but results in large model error. Models selected with smaller $\lambda$ values tend to result in lower model error and in turn appear smoother. This trade-off highlights a strength of our framework as it offers interpretable options, enabling experts to balance fidelity to data against model simplicity depending on the specific needs of the application.

\vspace{-2ex}
\begin{equation}
\label{Vanderpol-Library-insufficient}
N =
\begin{bmatrix}
\dot{x}_1 \\
\dot{x}_2 
\end{bmatrix} - 
\begin{bmatrix}
x_1 & x_2  & x_1 x_2 & x_1^2 & x_2^2 \\
x_1 & x_2  & x_1 x_2 & x_1^2 & x_2^2
\end{bmatrix}
\bigcirc \theta,
\end{equation}

\begin{table}[H]\label{tab:vdp_discovered_models}
\centering
{\scriptsize
\setlength{\tabcolsep}{3pt}
\renewcommand{\arraystretch}{1.2}
\begin{tabularx}{\textwidth}{|c|c|c|c|c|c|X|}
\hline
Rank & $\lambda$ & $R$ & Val. Err & Model Err & State Err & Model Discovered \\
\hline
1 & $10^{2}$ & $10^{-4}$ & $1.13 \times 10^{-1}$ & $3.75$ & $5.77 \times 10^{-2}$ &
$\begin{aligned}
\dot{x}_1 &= 1.063\,x_2 \\
\dot{x}_2 &= -0.993\,x_1
\end{aligned}$ \\
\hline
2 & $10^{0}$ & $10^{-4}$ & $3.75 \times 10^{-1}$ & $5.58 \times 10^{-2}$ & $3.39 \times 10^{-1}$ &
$\begin{aligned}
\dot{x}_1 &= 0.997\,x_2 \\
\dot{x}_2 &= -0.772\,x_1 - 0.297\,x_1x_2
\end{aligned}$ \\
\hline
3 & $10^{-1}$ & $10^{-4}$ & $4.80 \times 10^{-1}$ & $9.67 \times 10^{-4}$ & $4.42 \times 10^{-1}$ &
$\begin{aligned}
\dot{x}_1 &= 1.007\,x_2 \\
\dot{x}_2 &= -0.764\,x_1 - 0.071\,x_1^2 \, - 0.361\,x_1x_2 + 0.072\,x_2^2
\end{aligned}$ \\
\hline
\end{tabularx}
}
\caption{Discovered models at various $\lambda$ and $R$ values. Each balances model complexity with accuracy, offering candidates that users may select based on application-specific trade-offs.}
\vspace{-3ex}
\end{table}

\section{Levenberg-Marquardt} \label{Appx: LM_opt}

\begin{algorithm}
\footnotesize
\caption{Levenberg-Marquardt Optimization}
\label{alg:lm_optimization}
\begin{algorithmic}[1]
\Procedure{LevenbergMarquardt}{$f, \nabla f, \nabla^2 f, \mathbf{x}, \alpha, \tau, \epsilon$}
    \State Initialize $\alpha$, $\tau_1$, $\tau_2$ \Comment{e.g., $\alpha=1$, $\tau_1=0.25$, $\tau_2=2.0$}
    \While{$\|g\| > \epsilon$}
        \State $g \gets \nabla f(\mathbf{x}), \quad H_\alpha \gets \nabla^2 f(\mathbf{x}) + \alpha I$ \Comment{Gradient and modified Hessian}
        \State Solve $H_\alpha d = -g$
        \State $f_\text{old}, f_\text{test} \gets f(\mathbf{x}), f(\mathbf{x} + d)$
        \State $\text{predicted} \gets -g^T d - \frac{1}{2} d^T H d$ \Comment{Predicted decrease}
        \State $\rho \gets (f_\text{old} - f_\text{test}) / \text{predicted}$ \Comment{Agreement with model}
        \If{$f_\text{test} < f_\text{old}$}
            \State $\mathbf{x} \gets \mathbf{x} + d$ \Comment{Accept step}
            \If{$\rho < 0.25$}
                \State $\alpha \gets \min(\alpha \cdot \tau_2, \alpha_{\max})$ \Comment{Bad agreement, increase $\alpha$}
            \ElsIf{$\rho > 0.75$}
                \State $\alpha \gets \max(\alpha \cdot \tau_1, \alpha_{\min})$ \Comment{Good agreement, decrease $\alpha$}
            \EndIf
        \Else
            \State $\alpha \gets \min(\alpha \cdot \tau_2, \alpha_{\max})$ \Comment{Reject step, increase $\alpha$}
        \EndIf
    \EndWhile
\EndProcedure
\Return $\mathbf{x}$ 
\end{algorithmic}
\end{algorithm}

\end{document}